\documentclass[12pt]{article}
\usepackage[margin=2cm]{geometry}
\usepackage{ewLM,subfiles,titling}

\pdfoutput=1
\makeatletter
\def\writefile[2]{}

\input{IWCMainAux.aux}
\input{IWCSupplementAux.aux}

\begin{document}

	\title{Inference Without Compatibility}
	
	\author{Michael Law \and Ya\hspace{-.1em}'\hspace{-.1em}acov Ritov \thanks{Supported in part by NSF Grant DMS-1712962.}}
	\date{
		University of Michigan\\
		\today
	}
	
	\renewcommand\footnotemark{}
	
	\maketitle
	
	\begin{abstract}
		We consider hypotheses testing problems for three parameters in high-dimensional linear models with minimal sparsity assumptions of their type but without any compatibility conditions.  Under this framework, we construct the first $\sqrt{n}$-consistent estimators for low-dimensional coefficients, the signal strength, and the noise level.  We support our results using numerical simulations and provide comparisons with other estimators.
	\end{abstract}

	\section{Introduction}\label{sectionintroduction}
	
	In the past decade, there has been much interest in high-dimensional linear models, particularly following the work of \citeA{tibshirani1996}.  However, it was not until the past few years that there have been methods to construct confidence intervals and p-values for particular covariates in the model.  Consider a high-dimensional partially linear model
	\begin{align}\label{equationplmy}
		Y = X \beta + \mu + \epsilon,
	\end{align}
	with $X\in \R^{n\times q}$, and $Y,\mu,\epsilon \in \R^n$.  In addition, we also observe covariates $Z\in\R^{n\times p}$ such that $\mu \approx Z\gamma$ for some sparse vector $\gamma \in \R^{p}$ (see Section \ref{sectionnotation} for details).  Regarding the size of each matrix, we assume that $q < n$ is fixed but $p > n$ is high-dimensional.  Our goal is to construct a confidence region for the entire vector $\beta$.
	
	In recent years, there have been mainly two approaches to constructing confidence intervals in high-dimensional linear models.  There have been approaches such as \citeA{lee2016}, which construct conditional confidence intervals for $\beta$ given that $\beta$ was selected by a procedure, such as the lasso.  Simultaneously, there has been work to construct unconditional confidence intervals for $\beta$, where $X$ is the a priori selected covariate of interest, such as \citeA{javanmard2014}, \citeA{vandegeer2014}, and \citeA{zhang2014}; the latter is also our focus.  To avoid digressions, we will not elaborate on the former.  A review of many of the current methods is available in \citeA{dezeure2015}.  Much of the existing literature relies on using a version of the de-sparsified lasso introduced simultaneously by \citeA{javanmard2014}, \citeA{vandegeer2014}, and \citeA{zhang2014}.  The idea behind the existing approaches is to invert the KKT conditions of the lasso and perform nodewise lasso to approximate the inverse covariance matrix of the design, which attempts to correct the bias introduced by the lasso.
	
	Since the lasso forms the basis for the procedure, certain assumptions must be made in order to ensure that the lasso enjoys the nice theoretical properties that have been developed over the past two decades.  The paper by \citeA{vandegeer2009} provides an overview of various assumptions that have been used to prove oracle inequalities for the lasso.  These assumptions are a consequence of the fact the lasso is used rather than being needed for the statistical problem.  In particular, for confidence intervals, \citeA{vandegeer2014} assume that the compatibility condition holds for the Gram matrix, which is the weakest assumption from \citeA{vandegeer2009}, and is essentially a necessary assumption for the lasso to enjoy the fast rate (cf \citeA{bellec2018}).  To quote the popular book by \citeA{buhlmann2011}, ``In fact, a compatibility condition	is nothing else than simply an assumption that makes our proof go through.''  However, this raises an important question on necessity:  Is the compatibility condition necessary for constructing confidence intervals in high-dimensions?
	
	The main contribution of this paper is proving that the compatibility condition or any of its variants is indeed not necessary for the statistical problem.  To this end, we provide an estimator which does not require the compatibility condition but still attains the semi-parametric efficiency bound.  Our assumption regarding sparsity is slightly stronger than the minimax rate required by \citeA{javanmard2018} since we allow a broader class of designs.  In particular, we show that, in the absence of compatibility, the rate established by \citeA{javanmard2018} is not attainable and a stronger sparsity assumption is required.
	
	There is also the recent work of \citeA{chernozhukov2018}, who consider the general problem of conducting inference on low-dimensional parameters with high-dimensional nuisance parameters.  One application of their general theory is for high-dimensional partially linear models, which is also our problem of interest.  A further discussion of their procedure is given in Remark \ref{remarkchernozhukov} below.
	
	As a consequence of our estimation procedure for $\beta$, we are able to construct a $\sqrt{n}$-consistent estimator of the signal strength and the noise variance, which we denote by $\sigmamu$ and $\sigma_\epsilon^2$ respectively, also without the compatibility condition.  The paper by \citeA{reid2016} provides a nice overview of different proposals for estimation of $\sigma_\epsilon^2$ using the lasso.  An early work in this direction is \citeA{fan2012}, who construct asymptotic confidence intervals for $\sigma_\epsilon^2$ under a sure screening property of the covariates; in the setting of the lasso, this requires a $\beta$-min condition.  \citeA{dicker2014} consider a similar problem of variance estimation using moment estimators that do not require sparsity of the underlying signal.  However, they do not consider the ultra high-dimensional setting nor the problem of inference.  Later, \citeA{janson2017} considered inference on the signal-to-noise ratio but the theory developed only applies to Gaussian designs.  For the problem of inference for $\sigmamu$, the work most similar with ours is \citeA{cai2018}, who consider a more general problem in the semi-supervised setting, but their results for the supervised framework require minimal non-zero eigenvalues on the covariance matrix.  To this end, we construct estimators that attain asymptotic variances equal to that of the efficient estimator in low-dimensions.
	
	For both problems, our approach involves using exponential weighting to aggregate over all models of a particular size.  Prima facie, this is a computationally hard problem but can be well approximated in practice.  To this end, we propose an algorithm inspired by \citeA{rigollet2011}.
	
	\subsection{Organization of the Paper}
	We will end the current section with the notation that will be used throughout the paper.  In Section \ref{sectionbeta}, we discuss the problem of conducting inference for low-dimensional $\beta$ in the presence of a high-dimensional nuisance vector $\mu$.  The setting of univariate $\beta$ is considered separately in Section \ref{sectionq=1} to motivate the general multivariate procedure of Section \ref{sectionq>1}.  We take a slight detour in Section \ref{sectionq=1correlated} to consider inference when the errors are correlated.  The section ends with a discussion on the necessity of the sparsity assumption in Section \ref{sectionsparsity}.  Then, in Section \ref{sectionsigmamu} and Section \ref{sectionsigmaepsilon}, we consider the problems of inference for $\sigmamu$ and $\sigma_\epsilon^2$ respectively.  In Section \ref{sectionnumerics},  we provide an overview of the computation of the estimators, which we apply in Section \ref{sectionsimulations} for numerical simulations.  The proofs for Sections \ref{sectionq=1} and \ref{sectionsparsity} are provided in Section \ref{sectionproofs}.  Additional simulation tables and the proofs for the remaining results are available in the Supplement.

	\subsection{General Notation and Definitions}\label{sectionnotation}
	Throughout, all of our variables have a dependence on $n$, but when it should not cause confusion, this dependence will be suppressed.  For a general vector $a$ and a matrix $A$, $a_j$ will denote the $j$'th entry of $a$, $A_j$ the $j$'th column of $A$, and $A^{(j)}$ the $j$'th row of $A$.  Then, $\left\Vert a \right\Vert$ will denote the standard Euclidean norm, with the dimension of the space being implicit from the vector, $\left\Vert a \right\Vert_1$ the $L_1$-norm, and $\left\Vert a \right\Vert_0$ the $L_0$-norm.  Furthermore, $\left\Vert A \right\Vert$ will denote the operator norm and $\left\Vert A \right\Vert_\text{HS}$ the Hilbert-Schmidt norm.
	
	Before defining weak sparsity, we will need to introduce some notation.  For $u \in \N$, $\m_u$ will denote the collection of all models of $Z$ of size $u$.  That is,
	\begin{align*}
		\m_u \defined \left\{ m \subseteq \left\{ 1,\dots,p \right\} : |m| = u \right\}.
	\end{align*}
	Then, for each $m\in\m_u$, $Z_m$ will denote the $n\times u$ sub-matrix of $Z$ corresponding to the columns indexed by $m$.  Moreover, $P_m$ will denote the projection onto the column space of $Z_m$ and $\pp{m}$ the projection onto the orthogonal complement.  We can now state the definition of weak sparsity.
	\begin{definition}\label{definitionws}
		A sequence of vectors $\mu$ is said to satisfy the \emph{weak sparsity property relative to $Z$} with sparsity $s$ at rate $k$ if the set
		\begin{align*}
			\mathcal{S}_\mu \defined \left\{m\in\m_s : \left\Vert \pp{m}\mu \right\Vert^2 = o(k) \right\}
		\end{align*}
		is non-empty.  A set $S\in \mathcal{S}_\mu$ is said to be a \emph{weakly sparse set} for the vector $\mu$.
		
		If the sequence of vectors $\mu$ is random, then they satisfy the \emph{weak sparsity property relative to $Z$ in probability} with sparsity $s$ at rate $k$ if the set
		\begin{align*}
			\mathcal{S}_\mu = \left\{ m\in\m_s : \left\Vert \pp{m}\mu \right\Vert^2 = \op(k) \right\}
		\end{align*}
		is non-empty.  A set $S\in \mathcal{S}_\mu$ is said to be a \emph{weakly sparse set in probability} for the vector $\mu$.
	\end{definition}
	Finally, similar to other works on de-biased inference, we will consider sub-Gaussian errors, which is defined below.
	\begin{definition}\label{definitionsg}
		A mean zero random vector $\xi \in \R^n$ is said to be \emph{sub-Gaussian} with parameter $K$ if
		\begin{align*}
			\e \exp\left( \lambda^\T \xi \right) \leq \exp\left( \frac{ K^2 \left\Vert \lambda \right\Vert^2}{2} \right)
		\end{align*}
		for all vectors $\lambda \in \R^n$.
	\end{definition}

	\section{Inference for $\beta$}\label{sectionbeta}
	In this section, we consider the main problem of constructing confidence regions for $\beta$.  The model that we consider is given in equation \eqref{equationplmy}, which we reproduce below for convenience,
	\begin{align}\label{equationplmy2}
		Y = X\beta + \mu + \epsilon.
	\end{align}
	For this section, we will assume that $\mu$ satisfies the weak sparsity property relative to $Z$ at rate $\sqrt{n}$, but the results still hold if we assume the weak sparsity property in probability.  In addition to this partially linear model, we also assume that there exists matrices $N,H\in\R^{n\times q}$ such that each column of $X$ satisfies a partially linear model, denoted by
	\begin{align*}
		X_j = N_j + H_j,
	\end{align*}
	where $N_j$ satisfies the weak sparsity property relative to $Z$ at rate $\sqrt{n}$ for each $1\leq j \leq q$.  The weakly sparse set for each $N_j$ may be different, but the sparsity rate is uniformly $\sqrt{n}$.  In matrix form, we have that
	\begin{align}\label{equationplmx}
		X = N + H.
	\end{align}
We will assume that $H$ is sub-Gaussian. This assumption is certainly valid when $(X,Z)$ is jointly Gaussian. It precludes the case where $X$ is finitely supported, but is reasonable whenever the unbiased estimator of the risk used below is a good enough approximation for the risk of any given sub-model.
	By direct substitution, it follows that
	\begin{align*}
		Y = N\beta + \mu + H\beta + \epsilon.
	\end{align*}
	Then, since $\mu$ and each $N_j$ satisfy the weak sparsity property relative to $Z$ at rate $\sqrt{n}$, the vector $N\beta + \mu$ also satisfies the weak sparsity property relative to $Z$ at rate $\sqrt{n}$.
	
	In the case where $q=1$, we will write $\nu \defined N$ and $\eta \defined H$.  Moreover, we will define $\sigma_\epsilon^2 \defined \var(\epsilon_1)$, $\sigma_\eta^2 \defined \var(\eta_1)$ when $q=1$, and $\Sigma_H \defined \var\left(H^{(1)}\right)$ when $q>1$.
	
	\subsection{The Special Case:  $q=1$}\label{sectionq=1}
	Suppose that $q=1$ and let $S_\gamma$ and $S_\delta$ be weakly sparse sets for $\mu$ and $\nu$ respectively.  To motivate our procedure, we will assume temporarily that the models are in fact low-dimensional linear models, the set $S \defined S_\gamma \cup S_\delta$ is known, and $\epsilon \sim \n_n \left( 0_n , \sigma_\epsilon^2 I_n \right)$.  In particular, we are considering the low-dimensional linear models
	\begin{align*}
		Y &= X\beta + Z_{S_\gamma}\gamma + \epsilon = Z_S\theta + \eta\beta + \epsilon,\\
		X &= Z_{S_\delta}\delta + \eta,
	\end{align*}
	where $\theta = \delta\beta + \gamma$.  Then, by the Gauss-Markov Theorem, it is known that the efficient estimator in this low-dimensional problem is given by least-squares, which may be framed as the following three stage procedure:
	\begin{enumerate}
		\item Regress $Y$ on $Z_S$ using least-squares to obtain the fitted values $\hat{Y}$.
		\item Regress $X$ on $Z_S$ using least-squares to obtain the fitted values $\hat{X}$.
		\item Regress the residuals $Y - \hat{Y}$ on the the residuals $X - \hat{X}$ using least-squares to obtain the least-squares estimator $\betahatls$.
	\end{enumerate}
	In the high-dimensional setting, the first two stages can no longer be achieved using the classical least-squares approach.  However, since we are only interested in the fitted values $\hat{Y}$ and $\hat{X}$, this suggests using a high-dimensional prediction procedure to obtain the fitted values, and then applying low-dimensional least-squares on the residuals in the third stage.  The high-dimensional procedure that we will adopt is the exponential weights of \citeA{leung2006}, which has the salient feature of prediction consistency under very mild assumptions on the design.
	
	Before defining our estimators, we will state all of our assumptions.
	
	\begin{enumerate}[label=(A\arabic*)]
		\item \label{assumptionq=1boundednorms} The means $\mu$ and $\nu$ have squared norms that are $\Op(n)$.
		
		\item \label{assumptionq=1errordist} The entries of $\eta$ and $\epsilon$ are mutually independent and also independent of $Z$.  Moreover, the entries of $\eta$ and $\epsilon$ are each identically distributed sub-Gaussians with parameters $K_\eta$ and $K_\epsilon$ respectively.
		
		\item \label{assumptionq=1sparsity} The means $\mu$, $\nu$, and $\nu\beta + \mu$ are weakly sparse relative to $Z$ with sparsities $s_\gamma$, $s_\delta$, and $s_\theta$ respectively at rate $\sqrt{n}$.  Furthermore, the chosen sequence of sparsities satisfy $u_\gamma\geq s_\gamma$, $u_\delta \geq s_\delta$, and $u_\theta \geq s_\theta$ for $n$ sufficiently large and $\max\left(u_\gamma,u_\delta,u_\theta\right) = o(\sqrt{n}/\log(p))$.
	\end{enumerate}
	
	Now, we may define two sets of exponential weights, $w_{m,Y}$ and $w_{m,X}$, to estimate $\hat{Y}$ and $\hat{X}$ respectively.  Let
	\begin{align*}
		w_{m,Y} \defined \frac{\exp\left(-\frac{1}{\alpha_Y} \left\Vert \pp{m} Y \right\Vert^2\right)}{\sum_{k\in\m_{u_\theta}} \exp\left( -\frac{1}{\alpha_Y} \left\Vert \pp{k} Y \right\Vert^2 \right)}
	\end{align*}
	with $\alpha_Y > 4K_\epsilon^2$.
	
	\begin{remark}
		The exponential weights defined above do not subtract off the rank of the projection in the exponent as in \citeA{leung2006} since we only consider models of size $u_\theta$; the rank will cancel from the numerator and the denominator.
	\end{remark}
	
	Now, let $\thetahatm \defined \left( Z_m^\T Z_m \right)^{-1} Z_m^\T Y$ be the least-squares estimator for $\theta$ using the covariates $Z_m$.  Here, the matrix inverse is to be interpreted in a generalized sense if $Z_m^\T Z_m$ is rank deficient.  We will identify $\thetahatm$ with a vector in $\R^p$, with the support of $\thetahatm$ being indexed by $m$.  Then, we may estimate $\theta$ by
	\begin{align*}
		\thetahatew \defined \sum_{m\in\m_{u_\theta}} w_{m,Y} \thetahatm,
	\end{align*}
	with the prediction $\hat{Y}$ given by
	\begin{align*}
		\hat{Y} = Z \thetahatew.
	\end{align*}
	Similarly, we will define
	\begin{align*}
		w_{m,X} \defined \frac{\exp\left(-\frac{1}{\alpha_X} \left\Vert \pp{m} X \right\Vert^2\right)}{\sum_{k\in\m_{u_\delta}} \exp\left( -\frac{1}{\alpha_X} \left\Vert \pp{k} X \right\Vert^2 \right)},
	\end{align*}
	with $\alpha_X > 4K_\eta^2$.  Letting $\deltahatm$ denote the least-squares estimator of $\delta$ using the covariates $Z_m$ and identifying it with a vector in $\R^p$, we may define
	\begin{align*}
		\deltahatew \defined \sum_{m\in\m_{u_\delta}} w_{m,X} \deltahatm.
	\end{align*}
	Then, the fitted values of $X$ will be
	\begin{align*}
		\hat{X} = Z\deltahatew.
	\end{align*}
	Finally, for the last stage, the regression of $Y - Z\thetahatew$ on $X - Z\deltahatew$ will be given by
	\begin{align*}
		\betahatew \defined \frac{\left( X - Z\deltahatew \right)^\T \left( Y - Z\thetahatew \right)}{\left\Vert X - Z\deltahatew \right\Vert^2}.
	\end{align*}
	
	Before stating our main result, we will state a proposition regarding exponential weighting with sub-Gaussian errors.
	\begin{proposition}\label{propositionsubgaussianoracle}
		Consider a high-dimensional linear model given by
		\begin{align*}
			Y = \mu + \xi,
		\end{align*}
		for $\xi$ sub-Gaussian with parameter $K$.  Assume that $\mu$ is weakly sparse relative to $Z$ with sparsity $s$ and that $\limsup_{n\to\infty} \left\Vert \mu \right\Vert^2 = \mathcal{O}(n)$.  Assume further that the chosen sequence of sparsities $u\geq s$ satisfy $u = o(n^\tau/\log(p))$.  Letting $\hat{\gamma}_m$ denote the least-squares estimator for $\gamma$ using the covariates $Z_m$, define the exponential weights as
		\begin{align*}
			w_{m} \defined \frac{\exp\left( -\frac{1}{\alpha} \left\Vert \pp{m} Y \right\Vert^2 \right)}{\sum_{k\in\m_u} \exp\left( -\frac{1}{\alpha} \left\Vert \pp{k} Y \right\Vert^2 \right)},
		\end{align*}
		with $\alpha > 4K_\xi^2$.  Then,
		\begin{align*}
			\e \left\Vert \sum_{m\in\m_u} w_m Z\hat{\gamma}_m - \mu \right\Vert^2 = o(n^\tau).
		\end{align*}
	\end{proposition}
	\begin{remark}
		We would like to remark that the choice of $\alpha$ is consistent with \citeA{leung2006}.  In particular, when $\xi \sim \n_n\left(0_n,\sigma_\xi^2 I_n\right)$, the sub-Gaussian parameter is $K^2 = \sigma_\xi^2$, which gives the requirement that $\alpha > 4\sigma_\xi^2$.  In this setting, we would like to emphasize that the required value of $\alpha$ defies a Bayesian interpretation since the Bayes procedure requires a leading constant of $2$, as shown by \citeA{leung2006}.
	\end{remark}
	
	\begin{remark}
		The assumption that $\limsup_{n\to\infty} \left\Vert \mu \right\Vert^2 = \mathcal{O}(n)$ can be relaxed to hold in probability by weakening the conclusion to hold in probability rather than expectation (cf Corollary \ref{corollarybiasother}).
	\end{remark}
	
	For the remainder of the paper, we will only consider the setting where $\tau=1/2$.  As an immediate corollary, we have the following.
	\begin{corollary}\label{corollarysubgaussianoracle}
		Consider the models given in equations (\ref{equationplmy2}) and (\ref{equationplmx}) with $q=1$.  Under assumptions \ref{assumptionq=1boundednorms} -- \ref{assumptionq=1sparsity},
		\begin{align*}
			&\left\Vert \nu\beta + \mu - Z\thetahatew \right\Vert^2 = \op(\sqrt{n}),\\
			&\left\Vert \nu - Z\deltahatew \right\Vert^2 = \op(\sqrt{n}).
		\end{align*}
	\end{corollary}
	
	Finally, we can state the main result for $\betahatew$.
	\begin{theorem}\label{theoremq=1asympdist}
		Consider the models given in equations \eqref{equationplmy2} and \eqref{equationplmx} with $q=1$.  Under assumptions \ref{assumptionq=1boundednorms} -- \ref{assumptionq=1sparsity},
		\begin{align*}
			\sqrt{n}\left( \betahatew - \beta \right) \cond \n\left( 0 , \frac{\sigma_\epsilon^2}{\sigma_\eta^2} \right).
		\end{align*}
	\end{theorem}
	We would like to note that $\betahatew$ attains the information bound for estimating $\beta$ (cf Example 2.4.5 of \citeA{bickel1993} and Section 2.3.3 of \citeA{vandegeer2014}).
	
	\begin{remark}\label{remarkchernozhukov}
		The estimator, $\betahatew$, at first glance seems similar to the double/de-biased estimator of \citeA{chernozhukov2018} by considering exponential weighting as the estimation procedure for the propensity function.  However, the primary difference is that we do not rely on cross fitting to estimate the conditional mean of $X$ and $Y$ given the covariates $Z$.  Therefore, $\betahatew$ does not fall within the general framework of \citeA{chernozhukov2018} since exponential weighting intrinsically solves the high-dimensional in-sample prediction problem as opposed to the out-of-sample prediction problem.  That is, there is no guarantee that, after sample splitting, the out-of-sample predictions for the mean vector are necessarily consistent.
	\end{remark}
	
	To construct confidence intervals, we will need to estimate both $\sigma_\epsilon^2$ and $\sigma_\eta^2$.  We will defer explicitly defining estimators for the variance until Section \ref{sectionsigmaepsilon} but let $\sigmaepsilonhat$ and $\sigmaetahat$ be any of the three estimators proposed by Theorem \ref{theoremsigmaepsilon} for estimating variance.  Then, an asymptotic $(1-\alpha)$ confidence interval for $\beta$ is given by
	\begin{align*}
		\left( \betahatew - z_{\alpha/2} \sqrt{\frac{\sigmaepsilonhat}{\sigmaetahat n}} , \betahatew + z_{\alpha/2} \sqrt{\frac{\sigmaepsilonhat}{\sigmaetahat n}} \right),
	\end{align*}
	where $z_{\alpha/2}$ denotes the $\alpha/2$ upper quantile of the standard Gaussian distribution.

	\subsection{Correlated Gaussian Errors}\label{sectionq=1correlated}
	In this section, we take a slight detour away from classical high-dimensional partially linear models and consider the setting where the errors, $\epsilon$, are Gaussian but not necessarily independent and identically distributed.  The goal is to conduct inference on $\beta$, but, for simplicity, we will only consider the setting where $q=1$.  This model arises naturally if the model was a linear mixed model given by
	\begin{align*}
		Y = X\beta + \mu + W\zeta + \xi,
	\end{align*}
	where $\zeta$ are Gaussian random effects and $\xi$ is independent Gaussian noise.  \citeA{bradic2019} and \citeA{li2019} consider more general problems of testing fixed effects in high-dimensional linear mixed models, whereas we simply view the problem as a linear model with correlated noise.  Even when the errors are correlated, $\betahatew$ still has a Gaussian limit under proper rescaling.  Before stating the theorem, we will slightly modify assumption \ref{assumptionq=1errordist} to the setting where $\epsilon$ is correlated:
	\begin{enumerate}[label=(A\arabic*\,${}^\ast$)]
		\setcounter{enumi}{1}
		\item \label{assumptionq=1errordistcorrelated} The entries of $\eta \sim \n_n \left( 0 , \sigma_\eta^2 I_n \right)$ are independent of $Z$ and $\epsilon$.  The vector $\epsilon \sim \n_n \left( 0 , \Sigma_\epsilon \right)$ is independent of $Z$ with $\left\Vert \Sigma_\epsilon \right\Vert = \mathcal{O}(1)$ and
		$\Tr(\Sigma_\epsilon)/n \to \bar{d} > 0$.
	\end{enumerate}
	
	Now, we may state the main result for $\betahatew$ under correlation.
	
	\begin{theorem}\label{theoremq=1asympdistcorrelated}
		Consider the models given in equations \eqref{equationplmy2} and \eqref{equationplmx} with $q=1$.  Under assumptions \ref{assumptionq=1boundednorms}, \ref{assumptionq=1errordistcorrelated}, and \ref{assumptionq=1sparsity},
		\begin{align*}
			\sqrt{n}\left( \betahatew - \beta \right) \cond \n\left( 0 , \frac{\bar{d}}{\sigma_\eta^2} \right).
		\end{align*}
	\end{theorem}
	
	Again, we will defer defining an estimator for $\bar{d}$ and $\sigma_\eta^2$ until Section \ref{sectionsigmaepsilon}, in particular Corollary \ref{corollarysigmaepsiloncorrelated}.  Similar to the previous section, we may now construct confidence intervals for $\beta$ under this setting of correlation.

	\subsection{The General Case:  $q>1$}\label{sectionq>1}
	In the general setting where $q>1$, we may still rely on the perspective of high-dimensional prediction.  In particular, for $1 \leq j \leq q$, we may let $\deltahatewj$ denote the analogue of $\deltahatew$ for regressing $X_j$ on $Z$ and estimate $X_j$ by $Z\deltahatewj$.  Let $\Deltahatew \in \R^{p\times q}$ denote the matrix with columns given by $\deltahatewj$ for $1\leq j \leq q$.  Then, the multidimensional analogue of $\betahatew$ from Section \ref{sectionq=1} is given by
	\begin{align*}
		\betahatew \defined \left(\left( X - Z\Deltahatew \right)^\T \left( X - Z\Deltahatew \right)\right)^{-1} \left( X - Z\Deltahatew \right)^\T \left( Y - Z\thetahatew \right).
	\end{align*}
	We would like to emphasize that the definition here is identical to that given in Section \ref{sectionq=1} when $q=1$.
	
	Then, we will make the following assumptions.
	\begin{enumerate}[label=(B\arabic*)]
		\item \label{assumptionq>1boundednorms} The mean vectors $\mu$ and $N_j$ for $1\leq j \leq q$ have squared norms that are uniformly $\Op(n)$.
		
		\item \label{assumptionq>1errordist} The rows of $H$ and the entries of $\epsilon$ are independent and also independent of $Z$.  Moreover, the entries of the rows of $H$ and the entries of $\epsilon$ are each identically distributed sub-Gaussian with parameters $K_{\eta,j}$ and $K_\epsilon$ respectively.  Furthermore, $\Sigma_H$ is an invertible matrix.
		
		\item \label{assumptionq>1sparsity} All the mean vectors $\mu$, $N_j$ for $1\leq j \leq q$, and $N\beta + \mu$ are weakly sparse relative to $Z$ with sparsities $s_\gamma$, $s_{\delta,j}$ for $1\leq j \leq q$, and $s_\theta$ respectively at rate $\sqrt{n}$.  Furthermore, the chosen sequence of sparsities satisfy $u_\gamma\geq s_\gamma$, $u_{\delta,} \geq s_{\delta,j}$ for $1 \leq j \leq q$, and $u_\theta \geq s_\theta$ for $n$ sufficiently large and $\max\left(u_\gamma,\max_{1\leq j \leq q}\left(u_{\delta,j}\right),u_\theta\right) = o(\sqrt{n}/\log(p))$.
	\end{enumerate}
	
	We can now state the asymptotic distribution for $\betahatew$.
	\begin{theorem}\label{theoremq>1asympdist}
		Consider the models given in equations \eqref{equationplmy2} and \eqref{equationplmx}.  Under assumptions \ref{assumptionq>1boundednorms} -- \ref{assumptionq>1sparsity},
		\begin{align*}
			\sqrt{n}\left( \betahatew - \beta \right) \cond \n_q\left( 0_q , \sigma_\epsilon^2 \Sigma_H^{-1} \right).
		\end{align*}
	\end{theorem}
	
	Similar to before, to construct confidence regions, we will need to estimate $\Sigma_H$.  Therefore, we will consider
	\begin{align*}
		\SigmaHhat \defined \frac{1}{n}\left( X - Z\Deltahatew \right)^\T \left( X - Z\Deltahatew \right).
	\end{align*}
	This leads to the following proposition.
	
	\begin{proposition}\label{propositionq>1varianceestimates}
		Consider the models given in equations \eqref{equationplmy2} and \eqref{equationplmx}.  Under assumptions \ref{assumptionq>1boundednorms}, \ref{assumptionq>1errordist}, and \ref{assumptionq>1sparsity},
		\begin{align*}
			\SigmaHhat \conp \Sigma_H.
		\end{align*}
	\end{proposition}
	
	Then, an asymptotic $(1 - \alpha)$ confidence region for $\beta$ is given by
	\begin{align*}
		\left\{ \beta \in \R^q : \frac{n}{\sigmaepsilonhat} \left( \betahatew - \beta \right)^\T \SigmaHhat \left( \betahatew - \beta \right) \leq \chi^2_{q,\alpha} \right\},
	\end{align*}
	where $\chi^2_{q,\alpha}$ denotes the $\alpha$ upper quantile of a $\chi^2_q$ random variable.

	\subsection{Necessity of Sparsity Assumption}\label{sectionsparsity}
	In Section \ref{sectionq=1}, it was assumed that both $\mu$ and $\nu$ are weakly sparse with sparsity $s_\gamma$ and $s_\delta$ respectively at rate $\sqrt{n}$ in order for $\betahatew$ to have an asymptotic Gaussian distribution.  For simplicity, in the ensuing discussion, we will only consider the case where $q=1$, that there exists an $S\in \mathcal{S}_\mu$ such that $\left\Vert \pp{S} \mu \right\Vert^2 = 0$, and the design $(X,Z)$ is fully Gaussian with population covariance matrix $\Sigma$.  That is, $\Sigma = \var(X_1,Z^{(1)})$.  We will write $\Sigma_{Z,Z}$ to denote the $p\times p$ sub-block of $\Sigma$ corresponding to $Z$.  Letting $\Omega = \Sigma^{-1}$, it follows that
	\begin{align*}
		s_{\delta} = \left| \left\{ 1\leq j \leq p : \Omega_{1,j} \neq 0 \right\} \right|,
	\end{align*}
	which is equivalent to $s_\Omega$ from \citeA{javanmard2018}.  Compared to the de-biased lasso, \citeA{javanmard2018} showed that, if $s_\gamma = o(n/\log^2(p))$ and $\min(s_\gamma,s_\delta) = o(\sqrt{n}/\log(p))$, then the de-biased lasso has an asymptotic Gaussian distribution.  However, $\betahatew$ is a valid estimator on a larger class of designs, in particular incompatible designs, and Theorem \ref{theoremjrconjecture} formalizes this trade-off between sparsity and compatibility.  Before stating the theorem, we will need to introduce a bit of notation regarding our parameter space $\Theta$, which is defined as
	\begin{align*}
		\Theta(s_\gamma,s_\delta) \defined \{& \vartheta = \left( \beta , \gamma, \delta, \Sigma_{Z,Z}, \sigma_\eta^2, \sigma_\epsilon^2 \right)
		:\left\Vert \gamma \right\Vert_0 \leq s_\gamma, \left\Vert \delta \right\Vert_0 \leq s_\delta ,\\&
		\max\left(\gamma^\T \Sigma_{Z,Z} \gamma, \delta^\T \Sigma_{Z,Z} \delta,\sigma_\eta^2,\sigma_\epsilon^2\right) = \mathcal{O}(1)\}.
	\end{align*}
	
	\begin{theorem}\label{theoremjrconjecture}
		For $\vartheta \in \Theta(s_\gamma,s_\delta)$, consider the following model
		\begin{align*}
			&Z^{(1)},\dots,Z^{(n)} \overset{i.i.d}{\sim} \n_p\left(0_p, \Sigma_{Z,Z} \right),\\
			&\epsilon \sim \n_n\left( 0_n , \sigma_\epsilon^2 I_n \right),\\
			&\eta \sim \n_n \left( 0_n , \sigma_\eta^2 I_n \right),\\
			&Y = X\beta + Z\gamma + \epsilon,\\
			&X = Z\delta + \eta.
		\end{align*}
		Assume that either $s_\gamma = o(\sqrt{n}/\log(p))$ or $s_\delta = o(\sqrt{n}/\log(p))$.  If there exists a $\sqrt{n}$-consistent estimator of $\beta$ for all $\vartheta\in\Theta(s_\gamma,s_\delta)$, then both $s_\gamma = \mathcal{O}(\sqrt{n}/\log(p))$ and $s_\delta = \mathcal{O}(\sqrt{n}/\log(p))$.
	\end{theorem}
	
	In light of the results of \citeA{javanmard2018}, to construct a $\sqrt{n}$-consistent estimator of $\beta$, it must be the case that either $s_\gamma = o(\sqrt{n}/\log(p))$ or $s_\delta = o(\sqrt{n}/\log(p))$.  The previous theorem implies that the other sparsity must satisfy $\mathcal{O}(\sqrt{n}/\log(p))$.  Assumption \ref{assumptionq=1sparsity} is only mildly stronger, requiring $\max\left( s_\gamma, s_\delta \right) = o(\sqrt{n}/\log(p))$.

	\section{Inference for $\sigmamu$ and $\sigma_\epsilon^2$}\label{sectionsigmamuepsilon}
	
	In this section, we consider the problem of conducting inference for both $\sigmamu$ and $\sigma_\epsilon^2$.  \citeA{dicker2014}, \citeA{janson2017}, and \citeA{cai2018} provide interesting applications of both estimation and inference to which we refer the interested reader.  The main model that we consider is slightly different than that considered in the previous section.  Since we are not interested in the contribution of any particular covariate, we do not need to distinguish $X$ from $Z$.  Hence, we will set $q=0$ and consider the following model,
	\begin{align}\label{equationplmy3}
		Y = \mu + \epsilon.
	\end{align}
	Unlike Section \ref{sectionbeta}, we view $\mu$ as a random quantity, with $\sigmamu \defined \var(\mu_1)$.  Thus, $\sigmamu$ can be viewed as the explained variation in the data using the covariates $Z$.  Throughout this section, $S_\gamma$ will denote the weakly sparse set for $\mu$ with sparsity $s_\gamma$.  When constructing a $\sqrt{n}$-consistent estimator for $\sigmamu$, the asymptotic distribution will depend on the variance of $\mu_1^2$, which we will denote by $\kappa_\mu \defined \var\left(\mu_1^2\right)$.  Similarly, we will need to let $\kappa_\epsilon \defined \var\left( \epsilon_1^2 \right)$ when constructing confidence intervals for $\sigma_\epsilon^2$.

	\subsection{Inference for $\sigmamu$}\label{sectionsigmamu}
	
	To motivate our high-dimensional procedure, we will start by considering the low-dimensional setting.  Letting $S_\gamma$ denote a weakly sparse set for $\mu$ relative to $Z$ and identifying $\gamma$ with a vector in $\R^{s_\gamma}$, we will temporarily consider the linear model
	\begin{align}\label{equationlmmu}
		Y = Z_{S_\gamma} \gamma + \epsilon.
	\end{align}
	The natural estimator for $\sigmamu$ is given by $n^{-1} \left\Vert P_{S_\gamma} Y \right\Vert^2$.  The following proposition shows that this natural estimator is in fact efficient for estimating $\sigmamu$ with Gaussian errors.
	
	\begin{proposition}\label{propositionsigmamuefficiency}
		Consider the model given in equation \eqref{equationlmmu}.  Assume that the design $Z_{S_\gamma}$ has full column rank and $s_\gamma < n$ is fixed.  Then, the estimator $n^{-1} \left\Vert P_{S_\gamma} Y \right\Vert^2$ is efficient for estimating $\sigmamu$.
	\end{proposition}
	
	From the Central Limit Theorem, it is immediate that
	\begin{align*}
		\sqrt{n} \left( n^{-1} \left\Vert P_{S_\gamma} Y \right\Vert^2 - \sigmamu \right) \cond \n\left( 0 , \kappa_\mu + 4\sigmamu \sigma_\epsilon^2 \right).
	\end{align*}
	
	In the high-dimensional setting, there are three natural extensions of this low-dimensional efficient estimator using exponential weighting.  The first idea is to view $P_{S_\gamma} Y$ as the predicted values of $Y$ and directly use take the squared norm of the predicted values given by exponential weighting.  For $m\in\m_{u_\gamma}$, let $\gammahatm$ denote the least-squares estimator for $\gamma$ using the covariates $Z_m$ and set
	\begin{align*}
		\muhat \defined \sum_{m\in\m_{u_\gamma}} w_{m,Y} Z_m \gammahatm,
	\end{align*}
	where $w_{m,Y}$ is defined in Section \ref{sectionq=1}.  Then, we may consider the estimator
	\begin{align*}
		\sigmamuhats \defined \frac{1}{n}\left\Vert \muhat \right\Vert^2.
	\end{align*}
	Alternatively, we may take the perspective that exponential weights concentrate well around the models with high predictive capacity, which would suggest aggregating the squared norms,
	\begin{align*}
		\sigmamuhatm \defined \frac{1}{n}\sum_{m\in\m_{u_\gamma}} w_{m,Y} \left\Vert P_m Y \right\Vert^2.
	\end{align*}
	The last estimator that we consider is inspired by the low-dimensional maximum likelihood estimator for $\sigma_\epsilon^2$ and the fact that $\var(Y_1) = \sigmamu + \sigma_\epsilon^2$:
	\begin{align*}
		\sigmamuhatl \defined \frac{1}{n} \left( \left\Vert Y \right\Vert^2 - \left\Vert Y - \muhat \right\Vert^2\right).
	\end{align*}
	Before stating the main results for these estimators, we will first provide all of our assumptions.
	\begin{enumerate}[label=(C\arabic*)]
		\item \label{assumptionsigmamuboundednorms} The mean vector $\mu$ has independent and identically distributed entries with finite fourth moment.

		\item \label{assumptionsigmamuerrordist} The entries of $\epsilon$ are independent of $Z$.  Moreover, the entries of $\epsilon$ are independent and identically distributed sub-Gaussians with parameter $K_\epsilon$.
		
		\item \label{assumptionsigmamusparsity} The vector $\mu$ is weakly sparse relative to $Z$ with sparsity $s_\gamma$.  Furthermore, the chosen sparsity $u_\gamma$ satisfies $u_\gamma = o(\sqrt{n}/\log(p))$ and $u_\gamma \geq s_\gamma$ for $n$ sufficiently large.
	\end{enumerate}
	
	Assumption \ref{assumptionsigmamuboundednorms} implies that $\left\Vert \mu \right\Vert^2 = \Op(n)$.  By Jensen's inequality, it is immediate that $\sigmamuhats \leq \sigmamuhatm \leq \sigmamuhatl$.  However, it turns out that, under the above assumptions, these estimators are asymptotically equivalent at the $\sqrt{n}$-rate.  Recall that $\kappa_\mu \defined \var(\mu_1^2)$.  The following theorem provides the asymptotic distribution of the three estimators.
	
	\begin{theorem}\label{theoremsigmamu}
		Consider the model given in equation \eqref{equationplmy3}.  Suppose that $\sigmamu > 0$.  Under assumptions \ref{assumptionsigmamuboundednorms} -- \ref{assumptionsigmamusparsity},
		\begin{align*}
			\sqrt{n} \left( \sigmamuhat - \sigmamu \right) \cond \n\left( 0 , \kappa_\mu + 4\sigma_\epsilon^2 \sigmamu \right).
		\end{align*}
		where $\sigmamuhat$ is either $\sigmamuhats$, $\sigmamuhatm$, or $\sigmamuhatl$.
	\end{theorem}
	
	Since our interest is mainly asymptotic, we will write $\sigmamuhat$ to denote generically one of the estimators for $\sigmamu$.  To construct confidence intervals for $\sigmamu$, we will need to estimate $\kappa_\mu$, which may be accomplished by considering
	\begin{align*}
		\kappamuhat \defined \frac{1}{n} \sum_{j=1}^n \left( \muhatj^2 - \sigmamuhat \right)^2.
	\end{align*}
	
	The following proposition shows that $\kappamuhat$ is a consistent estimator for $\kappa_\mu$.
	
	\begin{proposition}\label{propositionkappamu}
		Consider the model given in equation \eqref{equationplmy3}.  Under assumptions \ref{assumptionsigmamuboundednorms} -- \ref{assumptionsigmamusparsity},
		\begin{align*}
			\kappamuhat \conp \kappa_\mu.
		\end{align*}
	\end{proposition}
	Therefore, an asymptotic $(1-\alpha)$ confidence interval for $\sigmamu$ is given by
	\begin{align}\label{cimu}
		\left( \sigmamuhat - z_{\alpha/2} \sqrt{\frac{\kappamuhat + 4 \sigmaepsilonhat\sigmamuhat}{n}} , \sigmamuhat + z_{\alpha/2} \sqrt{\frac{\kappamuhat + 4 \sigmaepsilonhat\sigmamuhat}{n}} \right).
	\end{align}

	\subsection{Inference for $\sigma_\epsilon^2$}\label{sectionsigmaepsilon}
	In this section, we are interested in constructing confidence intervals for $\sigma_\epsilon^2$.  In the low-dimensional setting with Gaussian errors, an estimator for $\sigma_\epsilon^2$ is given by maximum likelihood, which may be written as
	\begin{align*}
		\hat{\sigma}_{\epsilon,\text{ML}}^2 = \frac{1}{n} \left\Vert Y - P_{S_\gamma} Y \right\Vert^2.
	\end{align*}
	From classical parametric theory, $\hat{\sigma}_{\epsilon,\text{ML}}^2$ is an efficient estimator for $\sigma_\epsilon^2$ that achieves the information bound.  A natural extension in the high-dimensional setting is to view $P_{S_\gamma} Y$ as the predicted value and consider the estimator
	\begin{align*}
		\sigmaepsilonhats \defined \frac{1}{n} \left\Vert Y - \muhat \right\Vert^2,
	\end{align*}
	where $\muhat$ is defined in Section \ref{sectionsigmamu}.  Recalling that $\var\left( Y_1 \right) = \sigma_\mu^2 + \sigma_\epsilon^2$,
	we may consider two more estimators of $\sigma_\epsilon^2$ in light of the results of Section \ref{sectionsigmamu}, which are
	\begin{enumerate}
		\item $$ \sigmaepsilonhatm \defined \frac{1}{n} \left\Vert Y \right\Vert^2 - \sigmamuhatm.$$
		\item $$ \sigmaepsilonhatl \defined \frac{1}{n} \left\Vert Y \right\Vert^2 - \sigmamuhats.$$
	\end{enumerate}
	Again, by Jensen's inequality, it is immediate that $\sigmaepsilonhats \leq \sigmaepsilonhatm \leq \sigmaepsilonhatl$.  Similar to before, these three estimators are asymptotically equivalent at the $\sqrt{n}$-rate and the following theorem provides the asymptotic distribution for all three.
	
	\begin{theorem}\label{theoremsigmaepsilon}
		Consider the model given in equation \eqref{equationplmy3}.  Suppose that $\sigmamu > 0$.  Under assumptions \ref{assumptionsigmamuboundednorms} -- \ref{assumptionsigmamusparsity},
		\begin{align*}
			\sqrt{n} \left( \sigmaepsilonhat - \sigma_\epsilon^2 \right) \cond \n\left( 0 , \kappa_\epsilon \right).
		\end{align*}
		where $\sigmaepsilonhat$ is one of $\sigmaepsilonhats$, $\sigmaepsilonhatm$, or $\sigmaepsilonhatl$.
	\end{theorem}
	
	This gives us an immediate corollary to estimating $\bar{d}$ from Section \ref{sectionq=1correlated}, which requires the following assumption:
	\begin{enumerate}[label=(C$\arabic*^\ast$)]
		\setcounter{enumi}{1}
		\item \label{assumptionsigmamuerrordistcorrelated} The vector $\epsilon \sim \n_n \left( 0 , \Sigma_\epsilon \right)$ is independent of $Z$ with $\left\Vert \Sigma_\epsilon \right\Vert = \mathcal{O}(1)$ and
		$\Tr(\Sigma_\epsilon)/n \to \bar{d} > 0$.
	\end{enumerate}
	\begin{corollary}\label{corollarysigmaepsiloncorrelated}
		Consider the model given in equation \eqref{equationplmy3}.  Under assumptions \ref{assumptionsigmamuboundednorms}, \ref{assumptionsigmamuerrordistcorrelated}, and \ref{assumptionsigmamusparsity},
		\begin{align*}
			\sigmaepsilonhats \conp \bar{d}.
		\end{align*}
	\end{corollary}
	
	\begin{remark}
		Currently, in this section, we have assumed that $q=0$ but the theory for all three estimators of $\sigma_\epsilon^2$ are still valid when $q>0$.  In this setting, $X\beta + \mu$ is weakly sparse relative to $(X,Z)$ with sparsity $s_\gamma$ at rate $\sqrt{n}$.  Therefore, by using exponential weighting with the design $(X,Z)$, the above theorem implies that all three estimators are consistent for $\sigma_\epsilon^2$.
	\end{remark}
	
	\begin{remark}
		In practice, one may consider a version of the three estimators dividing by $n-u_\gamma$ instead of $n$, consistent with the low-dimensional unbiased mean squared error estimator.  Asymptotically, since $u_\gamma = o(\sqrt{n})$, they will have the same asymptotic distribution but seem to have better performance empirically in finite sample.
	\end{remark}
		
	Again, since $\sigmaepsilonhats$, $\sigmaepsilonhatm$, and $\sigmaepsilonhatl$ are asymptotically equivalent, we will write $\sigmaepsilonhat$ to denote a generically any of the three estimators.  To construct confidence intervals for $\sigma_\epsilon^2$, we will need to estimate $\kappa_\epsilon$.  The estimator that we propose is similar to $\kappamuhat$, namely we will defined $\kappaepsilonhat$ as
	\begin{align*}
		\kappaepsilonhat \defined \frac{1}{n} \sum_{j=1}^n \left( \left(y_j - \muhatj\right)^2 - \sigmaepsilonhat \right)^2.
	\end{align*}
	Analogous to Proposition \ref{propositionkappamu}, the following provides the consistency of $\kappaepsilonhat$.
	\begin{proposition}\label{propositionkappaepsilon}
		Consider the model given in equation \eqref{equationplmy3}.  Under assumptions \ref{assumptionsigmamuboundednorms} -- \ref{assumptionsigmamusparsity},
		\begin{align*}
			\kappaepsilonhat \conp \kappa_\epsilon.
		\end{align*}
	\end{proposition}
	Therefore, an asymptotic $(1-\alpha)$ confidence interval for $\sigma_\epsilon^2$ is given by
	\begin{align}\label{ciepsilon}
		\left( \sigmaepsilonhat - z_{\alpha/2} \sqrt{\frac{\kappaepsilonhat}{n}} , \sigmaepsilonhat + z_{\alpha/2} \sqrt{\frac{\kappaepsilonhat}{n}} \right).
	\end{align}

	\section{Implementation}\label{sectionnumerics}
	
	In this section, we describe a method to approximate all of the proposed estimators.  Since all of our estimators are based on exponential weighting, we will only detail the task of estimating $\thetahatew$, with the others being analogous.  Then, the goal of approximating $\thetahatew$ can be split into the following two tasks:
	
	\begin{enumerate}
		\item Determining the values of the tuning parameters $\alpha_Y$ and $u_\theta$.
		
		\item Aggregating over ${p \choose u_\theta}$ models.
	\end{enumerate}

	We will start with the second task.  Suppose temporarily that values of $\alpha_Y$ and $u_\theta$ have been selected.  To aggregate the models, we will follow the Metropolis Hastings scheme of \citeA{rigollet2011}.  Our approach slightly differs from theirs since we restrict our attention to $u_\theta$-sparse models whereas they consider models of varying sizes.
	
	Conditional on the data, the values of $\thetahatew$ and $\thetahatm$ for each $m\in\m_{u_\theta}$ are fixed.  We may view $\m_{u_\theta}$ as the vertices of the Johnson graph $J(p,u_\theta, u_\theta - 1)$ (cf \citeA{godsil2013}).  Then, for each $m\in\m_{u_\theta}$, by assigning weight $w_{m,Y}$ to vertex $m$, the target $\thetahatew$ may be viewed as the expectation of the fixed estimators $\thetahatm$ over the graph $J(p,u_\theta, u_\theta - 1)$, conditional on the observed data.  Hence, by taking a random walk over $J(p,u_\theta, u_\theta - 1)$, we may approximate $\thetahatew$.
	
	Before describing the algorithm, we need to introduce a bit of notation.  For any model $m\in\m_{u_\theta}$, we will let $\k_m$ denote the neighbors of $m$, which is given by
	\begin{align*}
		\k_m \defined \left\{k\in\m_{u_\theta} : \left| k \cap m \right|=u_\theta - 1\right\}.
	\end{align*}
	Moreover, we will write $RSS_m \defined \left\Vert \pp{m} Y \right\Vert^2$, the residual sum of squares.  Finally, let $T_0$ denote some burn-in time for the Markov chain and $T$ denote the number of samples from the Markov chain.  This will yield the following algorithm, which closely parallels \citeA{rigollet2011}.
	\vspace{2ex}
		
	\begin{algorithm}[H]
		\SetAlgoLined
		\KwResult{Approximates $\thetahatew$}
		Initialize a random point $m_0\in\m_u$ and compute $RSS_{m_0}$\;
		\For{$t=1,\dots,T$}
		{
			Uniformly select $k\in\k_{m_t}$ and compute $RSS_k$\;
			Generate a random variable $m_{t+1}$ by
			\begin{align*}
			m_{t+1} = \begin{cases}
			m_{t} &\text{ with probability } \exp\left(-\frac{1}{\alpha_Y} (RSS_k - RSS_{m_t})\right)\text{\;}\\
			k &\text{ with probability } 1 - \exp\left(-\frac{1}{\alpha_Y} (RSS_k - RSS_{m_t})\right)\text{\;}
			\end{cases}
			\end{align*}
			\If{$t>T_0$}
			{Compute $\hat{\theta}_{t+1} \gets (Z_{m_{t+1}}^\T Z_{m_{t+1}})^{-1} Z_{m_{t+1}}^\T Y$, embedded as a vector in $\R^p$\;}
		}
		\Return
		\begin{align*}
			\frac{1}{T}\sum_{t=T_0+1}^{T_0+T} \hat{\theta}_{t+1}\text{\;}
		\end{align*}
		\caption{Exponential weighting}
	\end{algorithm}
	\vspace{2ex}
	
	Then, analogous to Theorem 7.1 of \citeA{rigollet2011}, it will follow that
	\begin{align*}
		\lim_{T\to\infty} \frac{1}{T}\sum_{t=T_0+1}^{T_0+T} \hat{\theta}_{t+1} = \thetahatew && \p \text{ almost surely}.
	\end{align*}
	
	Now, for the first task, we may construct a grid of parameter points and use cross-validation to jointly tune the parameters using the above algorithm.  Since both $\alpha_Y$ and $u_\theta$ do not need to be known exactly, but need to be tuned to be larger than a threshold, the grid can be quite coarse to ease the computational burden.
	
	Computation in the ultrahigh-dimension is inherently difficult.  In view of \citeA{zhang2014lower}, there is no polynomial time algorithm that achieves the minimax rate for prediction without the restricted eigenvalue condition.  However, we do not know any algorithm that verifies the restricted eigenvalue condition in polynomial time (cf \citeA{raskutti2010}).  In this paper, we completely avoid assuming a condition like the restricted eigenvalue condition and therefore we cannot guarantee polynomial time convergence.  Yet, the algorithm behaves well in practice, as can be seen from the simulations in the following section.
	
	\section{Simulations}\label{sectionsimulations}
	We divide this section into two parts, corresponding to simulations for $\beta$ and simulations for variance components $\sigma_\mu^2$ and $\sigma_\epsilon^2$.  Additional simulation tables are included in the Supplement.
	
	\subsection{Simulations for $\beta$}
	
	For ease of comparison, our simulations will be similar to those given in \citeA{vandegeer2014}.  For the linear models
	\begin{align*}
		&Y = X\beta + \mu + \epsilon,\\
		&X_j = N_j + H_j,
	\end{align*}
	we will consider the setting where $n=100$ and $p=500$.  There are a few parameters with which we will experiment: $q$, $\beta$, the distribution of the design and errors, the sparsities, and the signal to noise ratio.  For each parameter pairing, we run $500$ simulations.  All confidence intervals will be constructed at the nominal $95\%$ level.
	
	Since the number of parameters of interest is fixed and low-dimensional, we will consider the settings where $q\in\left\{1,3\right\}$.  To assess both the coverage and the power, we will let $\beta$ be a vector in $\R^q$ with values in $\left\{0,1\right\}$.  To experiment with the robustness to the sub-Gaussianity assumption, we will use Gaussian, double exponential, and $t(3)$ distributions for the errors, all scaled to have mean zero and unit variance.  We will denote these distributions by z, e, and t respectively.  Therefore, $\sigma_\epsilon^2 = 1$ throughout this section.  The design will have the same distribution as the error, but with an equi-correlation covariance matrix.  That is, we consider the covariance matrix, $\Sigma(Z)$ to be
	\begin{align*}
		\Sigma(Z)_{i,j} =
		\begin{cases}
			1 &\text{ if } i=j\\
			\rho &\text{ if } i\neq j
		\end{cases}
	\end{align*}
	for $\rho \in \left\{0,0.8\right\}$.  When $q=3$, the covariance matrix for $H^{(1)}$, denoted by $\Sigma(H)$, will also be equi-correlation,
	\begin{align*}
		\Sigma(H) = \begin{cases}
		\sigma_\eta^2& \text{ if } i=j.\\
		0.5\sigma_\eta^2 &\text{ if } i\neq j,
		\end{cases}
	\end{align*}
	where $\sigma_\eta^2$ is chosen so that $\var(X_1)=1$.
	
	Similar to \citeA{vandegeer2014}, we will let the sparsity $s_\gamma\in\left\{3,15\right\}$, and, for simplicity, set $s_\delta=s_\gamma$.  We will set the signal to noise ratio of $\mu$ to $\epsilon$, which is given by $\sigma_\mu^2/\sigma_\epsilon^2$, to be $2$.  Since large values of the signal to noise ratio ($SNR$) of $N_j$ to $H_j$ correspond to highly correlated designs, we will also consider $SNR_X \defined \sigma_\nu^2/\sigma_\eta^2 \in \{2,1000\}$.
	
	For our simulations, we will say $\mu$ is weakly sparse relative to $Z$ with sparsity $s_\gamma$ at rate $\sqrt{n}$ if there exists an $s_\gamma$-sparse set $S$ and vector $\gamma_S$ such that $\var(\mu_1 - \left(Z_S\gamma_S\right)_1) \leq n^{-1/2}$.  In particular, we will consider vectors $\gamma$ of the form
	\begin{align*}
		\gamma_j \propto \pi(j)^{-\kappa} && j=1,\dots,p
	\end{align*}
	for some value $\kappa >0$ and permutation $\pi : \left\{1,\dots, p\right\} \to \left\{1,\dots, p\right\}$.  A similar approach is applied for $\Delta$.
	
	We will compare our estimators with a few other procedures:
	\begin{enumerate}
		\item (LS) Oracle least-squares that knows the true weakly sparse set $S_\gamma$.
		
		\item (DLA) De-biased lasso from \citeA{dezeure2015} as implemented in the R package {\tt hdi}.  We only apply this when $q=1$.
		
		\item (SILM) Simultaneous inference for high-dimensional linear models of \citeA{zhang2017} as implemented in the R package {\tt SILM}.
		
		\item (DML) Double/de-biased machine learning of \citeA{chernozhukov2018} with $4$ folds using the scaled lasso of \citeA{sun2012} as the estimation procedure as implemented in the R package {\tt scalreg}.  We only apply this when $q=1$.
		
		\item ($\text{EW}_{I}$), ($\text{EW}_{II}$), ($\text{EW}_{III}$) Exponential weights using $\sigmaepsilonhats$, $\sigmaepsilonhatm$, and $\sigmaepsilonhatl$ respectively.  We tune the parameters using cross-validation with $T_0 = 3000$ and $T=7000$.
	\end{enumerate}

	To evaluate the procedures, we use the following two measures
	\begin{enumerate}
		\item (AvgCov) Average coverage:  The percentage of time the true value of $\beta$ falls inside the confidence region.
		\item (AvgLen) Average length:  The average length of the confidence interval (only when $q=1$).
	\end{enumerate}

	The results are given in Table 1 and Tables \ref{tableq3beta0desz}--\ref{tableq3beta1dest}  from the Supplement.  In the $q=1$ setting with $SNR_X=2$, the coverage is comparable amongst all of the estimators.  However, the de-biased lasso and the SILM procedure are slightly preferable in this regime since the length of the intervals are slightly shorter.  When $\beta=0$, $SNR_X=1000$, and $\rho=0.8$, the coverage of the de-biased lasso is quite poor, with less than a $25\%$ coverage against a nominal rate of $95\%$.  The result should not be surprising since this corresponds to a setting of high correlation in the design, which weakens the compatibility condition.  The double/de-biased machine learning approach has strong nominal coverage in this regime (about $100\%$), but the length of the intervals are significantly longer than the other procedures (about four to five times longer than exponential weighting).  When $\beta=1$, $SNR_X=1000$, and $\rho=0.8$, we note that the SILM procedure no longer maintains nominal coverage.  The results remain the same when we consider $q=3$ and different distributions for the design and the errors.
	
\begin{table}[]
\centering
\caption{Simulations for $\beta$ with Gaussian design and errors when q=1 and $\beta=$0} 
\label{tableq1beta0desz}
\begin{tabular}{|l|l|rrrrrrrr|}
   \hline
 & $snr_X$ & 2 & 2 & 2 & 2 & 1000 & 1000 & 1000 & 1000 \\ 
   & $\rho$ & 0 & 0 & 0.8 & 0.8 & 0 & 0 & 0.8 & 0.8 \\ 
   & $s_\delta,s_\gamma$ & 3 & 15 & 3 & 15 & 3 & 15 & 3 & 15 \\ 
   \hline
 & LS & 0.946 & 0.880 & 0.946 & 0.958 & 0.942 & 0.908 & 0.938 & 0.930 \\ 
   & DLA & 0.958 & 0.884 & 0.976 & 0.978 & 0.954 & 0.870 & 0.218 & 0.170 \\ 
   & SILM & 0.970 & 0.872 & 0.962 & 0.970 & 0.958 & 0.812 & 0.900 & 0.902 \\ 
  AvgCov & DML & 0.966 & 0.850 & 0.956 & 0.946 & 0.982 & 0.844 & 1.000 & 1.000 \\ 
   & $\text{EW}_{I}$ & 0.956 & 0.868 & 0.956 & 0.962 & 0.960 & 0.828 & 0.954 & 0.968 \\ 
   & $\text{EW}_{II}$ & 0.978 & 0.912 & 0.976 & 0.980 & 0.972 & 0.898 & 0.966 & 0.984 \\ 
   & $\text{EW}_{III}$ & 0.984 & 0.938 & 0.984 & 0.994 & 0.980 & 0.936 & 0.980 & 0.994 \\ 
   \hline
 & LS & 0.427 & 0.462 & 0.589 & 0.684 & 0.430 & 0.467 & 0.919 & 1.440 \\ 
   & DLA & 0.493 & 0.532 & 0.689 & 0.700 & 0.530 & 0.547 & 0.544 & 0.501 \\ 
   & SILM & 0.529 & 0.559 & 0.670 & 0.697 & 0.623 & 0.609 & 0.666 & 0.646 \\ 
  AvgLen & DML & 0.650 & 0.634 & 0.694 & 0.692 & 1.510 & 0.881 & 10.600 & 11.100 \\ 
   & $\text{EW}_{I}$ & 0.623 & 0.636 & 0.700 & 0.716 & 1.060 & 0.774 & 1.910 & 1.830 \\ 
   & $\text{EW}_{II}$ & 0.690 & 0.710 & 0.768 & 0.797 & 1.170 & 0.868 & 2.100 & 2.040 \\ 
   & $\text{EW}_{III}$ & 0.749 & 0.776 & 0.830 & 0.871 & 1.280 & 0.951 & 2.270 & 2.240 \\ 
   \hline
\end{tabular}
\end{table}

	\subsection{Simulations for $\sigma_\mu^2$ and $\sigma_\epsilon^2$}
	
	In this section, we set $q=0$ and only consider the setting of strong sparsity (ie. $\mu = Z\gamma$ for some vector $\gamma \in \R^p$ satisfying $\left\Vert \gamma \right\Vert_0 = s_\gamma$).  This reduces the linear model to
	\begin{align*}
		Y = Z\gamma + \epsilon.
	\end{align*}
	We still consider the setting where $n=100$ and $p=500$.  The value of $\sigma_\mu^2 =2$ and $\sigma_\epsilon^2=1$ throughout these simulations.  The parameters with which we will experiment are the distributions of the design and errors and the sparsity.
	
	Again, we will consider Gaussian, double exponential, and $t(3)$ distributions for the design and the errors.  The design will have an equi-correlation structure with $\rho\in\left\{0,0.8\right\}$ and the sparsity will satisfy $s_\gamma \in\left\{3,15\right\}$.
	
	The vector of coefficients, $\gamma$, will have $s_\gamma$ components generated from uniform(-1,1) and $p-s_\gamma$ components that are zero.  The values will then be scaled such that $\sigma_\mu^2 = \gamma^\T \Sigma_Z \gamma = 2$.
	
	For estimation of $\sigma_\mu^2$, we will compare our results with an oracular estimator based on low-dimensional least-squares and the recent proposal of CHIVE.
	
	\begin{enumerate}
		\item (LS) Oracle least-squares that knows the true strongly sparse set $S_\gamma$ using equation \eqref{cimu}.
		
		\item (CHIVE) The calibrated inference for high-dimensional variance explained of \citeA{cai2018}.  We follow Algorithm 1 of the paper with $\tau_0^2\in \left\{0,2,4,6\right\}$.
		
		\item ($\text{EW}_{I}$), ($\text{EW}_{II}$), ($\text{EW}_{III}$) Exponential weighting using $\sigmamuhats$, $\sigmamuhatm$, and $\sigmamuhatl$ respectively.  We tune the parameters using cross-validation with $T_0 = 3000$ and $T=7000$.
	\end{enumerate}
	
	The results are presented in Table \ref{tablemusgamma3} and Table \ref{tablemusgamma15} from the Supplement.  We note that the coverage of the least-squares procedure is close to the nominal $95\%$ rate when $s_\gamma=3$ and the errors are either Gaussian or double exponential.  The coverage is significantly worse for the $t(3)$ design, which should not be surprising since the fourth moment is not defined for this distribution.  However, when $s_\gamma=15$, the coverage of least-squares falls, which establishes a reference for the problem difficulty, since Proposition \ref{propositionsigmamuefficiency} establishes the efficiency of least-squares in this problem.
	
	Amongst the exponential weighting estimators, when $s_\gamma=3$ and the errors are Gaussian or double exponential, the procedure based on $\sigmamuhats$ has the best performance and $\sigmamuhatl$ has the coverage when the errors are $t$ distributed.  For higher sparsity, no one estimators dominates the others; depending on our assumptions, any of the three estimators may be preferable.  Compared with CHIVE, the best exponential weighting procedure seems to be able to achieve comparable coverage with significantly shorter intervals, which can be seen across all of our simulation settings.
		
\begin{table}[H]
\centering
\caption{Simulations for $\sigma_{\mu}^2$ with $s_{\gamma}=$3} 
\label{tablemusgamma3}
\begin{tabular}{|l|l|rrrrrr|}
   \hline
 & Distribution & z & z & e & e & t & t \\ 
   & $\rho$ & 0 & 0.8 & 0 & 0.8 & 0 & 0.8 \\ 
   \hline
 & LS & 0.922 & 0.948 & 0.914 & 0.934 & 0.808 & 0.802 \\ 
   & $\text{CHIVE}_{0}$ & 0.698 & 0.532 & 0.690 & 0.604 & 0.554 & 0.526 \\ 
   & $\text{CHIVE}_{2}$ & 0.818 & 0.668 & 0.792 & 0.702 & 0.712 & 0.634 \\ 
  AvgCov & $\text{CHIVE}_{4}$ & 0.888 & 0.748 & 0.848 & 0.762 & 0.770 & 0.704 \\ 
   & $\text{CHIVE}_{6}$ & 0.890 & 0.772 & 0.898 & 0.790 & 0.860 & 0.746 \\ 
   & $\text{EW}_{I}$ & 0.852 & 0.850 & 0.854 & 0.862 & 0.780 & 0.778 \\ 
   & $\text{EW}_{II}$ & 0.804 & 0.772 & 0.820 & 0.838 & 0.812 & 0.828 \\ 
   & $\text{EW}_{III}$ & 0.708 & 0.644 & 0.744 & 0.762 & 0.820 & 0.866 \\ 
   \hline
 & LS & 1.520 & 1.510 & 1.800 & 1.950 & 2.430 & 2.950 \\ 
   & $\text{CHIVE}_{0}$ & 0.998 & 0.937 & 1.160 & 1.190 & 1.670 & 2.130 \\ 
   & $\text{CHIVE}_{2}$ & 1.520 & 1.560 & 1.650 & 1.740 & 2.150 & 2.640 \\ 
  AvgLen & $\text{CHIVE}_{4}$ & 1.890 & 1.970 & 2.010 & 2.120 & 2.500 & 2.980 \\ 
   & $\text{CHIVE}_{6}$ & 2.210 & 2.300 & 2.310 & 2.440 & 2.780 & 3.270 \\ 
   & $\text{EW}_{I}$ & 1.470 & 1.440 & 1.750 & 1.850 & 2.390 & 2.840 \\ 
   & $\text{EW}_{II}$ & 1.420 & 1.390 & 1.710 & 1.810 & 2.370 & 2.810 \\ 
   & $\text{EW}_{III}$ & 1.370 & 1.320 & 1.670 & 1.760 & 2.340 & 2.780 \\ 
   \hline
\end{tabular}
\end{table}

	For the estimation of $\sigma_\epsilon^2$, we will consider the oracular least-squares, the scaled lasso estimator, and the refitted cross-validation with Sure Independence Screening, along with our proposed procedures based on exponential weighting.
	
	\begin{enumerate}
		\item (LS) Oracle least-squares that knows the true strongly sparse set $S_\gamma$ using equation \eqref{ciepsilon}.
		
		\item (SL) Scaled lasso as implemented in the R package {\tt scalreg} with a confidence interval constructed using Theorem 2 of \citeA{sun2012}.
		
		\item (RCV-SIS) Refitted cross-validation of \citeA{fan2012} using the Sure Independence Screening of \citeA{fan2008} as implemented in the R package {\tt SIS} in the first stage.  The confidence interval is constructed using Theorem 2 of \citeA{fan2012}, with $\e\epsilon^4$ estimated by Proposition \ref{propositionkappaepsilon} of the present paper.
		
		\item ($\text{EW}_{I}$), ($\text{EW}_{II}$), ($\text{EW}_{III}$) Exponential weighting using $\sigmaepsilonhats$, $\sigmaepsilonhatm$, and $\sigmaepsilonhatl$ respectively.  We tune the parameters using cross-validation with $T_0 = 3000$ and $T=7000$.
	\end{enumerate}

	The results are given in Table \ref{tableepsilonsgamma3} and Table \ref{tableepsilonsgamma15} from the Supplement.  When the signal is very sparse, $s_\gamma=3$, and there is no correlation in the design, scaled lasso has better coverage than exponential weighting.  However, as the correlation increases to $\rho=0.8$, the confidence intervals constructed using $\sigmaepsilonhatm$ outperforms scaled lasso both in terms of coverage and average length.  When the model is less sparse, $\sigmaepsilonhats$ has comparable or better performance than scaled lasso.  The poor performance of refitted cross-validation with Sure Independence Screening in the $s_\gamma=15$ case should not come as a surprise since the signal to noise ratio is kept constant.  The task of sure screening $15$ active covariates out of $500$ with low signal strength from $50$ observations is very difficult.
	
\begin{table}[H]
\centering
\caption{Simulations for $\sigma_{\epsilon}^2$ with $s_{\gamma}=$3} 
\label{tableepsilonsgamma3}
\begin{tabular}{|l|l|rrrrrr|}
   \hline
 & Distribution & z & z & e & e & t & t \\ 
   & $\rho$ & 0 & 0.8 & 0 & 0.8 & 0 & 0.8 \\ 
   \hline
 & LS & 0.938 & 0.912 & 0.952 & 0.940 & 0.918 & 0.912 \\ 
   & SL & 1.000 & 0.730 & 0.998 & 0.730 & 0.994 & 0.756 \\ 
  AvgCov & RCV-SIS & 0.684 & 0.646 & 0.688 & 0.644 & 0.638 & 0.606 \\ 
   & $\text{EW}_{I}$ & 0.616 & 0.608 & 0.678 & 0.674 & 0.650 & 0.690 \\ 
   & $\text{EW}_{II}$ & 0.862 & 0.828 & 0.872 & 0.846 & 0.852 & 0.814 \\ 
   & $\text{EW}_{III}$ & 0.672 & 0.458 & 0.660 & 0.488 & 0.636 & 0.430 \\ 
   \hline
 & LS & 0.532 & 0.529 & 0.545 & 0.528 & 0.534 & 0.534 \\ 
   & SL & 0.599 & 0.670 & 0.602 & 0.665 & 0.602 & 0.659 \\ 
  AvgLen & RCV-SIS & 0.485 & 0.509 & 0.508 & 0.514 & 0.554 & 0.539 \\ 
   & $\text{EW}_{I}$ & 0.430 & 0.427 & 0.442 & 0.438 & 0.435 & 0.447 \\ 
   & $\text{EW}_{II}$ & 0.441 & 0.444 & 0.453 & 0.453 & 0.446 & 0.463 \\ 
   & $\text{EW}_{III}$ & 0.462 & 0.475 & 0.473 & 0.480 & 0.466 & 0.492 \\ 
   \hline
\end{tabular}
\end{table}

	\section{Proofs}\label{sectionproofs}
	
	\subsection{Proofs for Section \ref{sectionq=1}}
	For ease of reference in later proofs, we will prove Proposition \ref{propositionsubgaussianoracle} as two lemmata.
	\begin{lemma}\label{lemmabiashw}
		Let $w_m$ be any collection of convex weights over $\m_u$, and $\xi$ be a sub-Gaussian vector with parameter $K_\xi$, independent of $Z$.  If $u = o(n^\tau/\log(p))$, then
		\begin{align*}
			\e \left( \sum_{m\in\m_u} w_m \left\Vert P_m \xi \right\Vert^2 \right) = o(n^\tau).
		\end{align*}
	\end{lemma}
	\begin{proof}
		Fix $t>0$ arbitrarily.  Define the event $\Tset_t$ as
		\begin{align*}
			\Tset_t \defined \bigcap_{m\in\m_u} \left\{ \left\Vert P_m \xi \right\Vert^2 \leq K_\xi^2 \left( u + 2\sqrt{utn^{\tau}} + 2tn^{\tau} \right) \right\}.
		\end{align*}
		For any fixed $m\in\m_u$, it follows from Theorem 2.1 of \citeA{hsu2012} that
		\begin{align*}
			\p\left( \left\Vert P_m \xi \right\Vert^2 > K_\xi^2 \left( u + 2\sqrt{utn^{\tau}} + 2tn^{\tau} \right) \right) \leq \exp\left( -tn^{\tau} \right).
		\end{align*}
		Therefore,
		\begin{align}\label{equationtsettcprob}
			\p \left( \Tset_t^\C \right) \leq \exp\left( -tn^\tau + \log\left( |\m_u|\right) \right).
		\end{align}
		Now, note that
		\begin{align*}
			\e \left( \sum_{m\in\m_u} w_m \left\Vert P_m \xi \right\Vert^2 \right)
			= \e \left( \sum_{m\in\m_u} w_m \left\Vert P_m \xi \right\Vert^2 \indic{\Tset_t} \right)
			+ \e \left( \sum_{m\in\m_u} w_m \left\Vert P_m \xi \right\Vert^2 \indic{\Tset_t^\C} \right).
		\end{align*}
		For the first term, by the definition of $\Tset_t$,
		\begin{align*}
			\limsup_{n\to\infty} n^{-\tau} \e \left( \sum_{m\in\m_u} w_m \left\Vert P_m \xi \right\Vert^2 \indic{\Tset_t} \right) \leq 2tK_\xi^2.
		\end{align*}
		For the second term, by Cauchy-Schwarz and equation \eqref{equationtsettcprob}, it follows that
		\begin{align*}
			\limsup_{n\to\infty}n^{-\tau}\e \left( \sum_{m\in\m_u} w_m \left\Vert P_m \xi \right\Vert^2 \indic{\Tset_t^\C} \right)
			&\leq \limsup_{n\to\infty}n^{-\tau} \e \left( \left\Vert \xi \right\Vert^2 \indic{\Tset_t^\C} \right)\\
			&\leq \limsup_{n\to\infty}n^{-\tau} \e \left( \left\Vert \xi \right\Vert^4 \right)^{1/2} \p\left( \Tset_t^\C \right)^{1/2}\\
			&=  0.
		\end{align*}
		Therefore,
		\begin{align*}
			\limsup_{n\to\infty} n^{-\tau} \e \left( \sum_{m\in\m_u} w_m \left\Vert P_m \xi \right\Vert^2 \right) \leq 2tK_\xi^2.
		\end{align*}
		Since $t>0$ was arbitrary, this finishes the proof.
	\end{proof}
	
	\begin{lemma}\label{lemmabiasother}
		Under the assumptions and setup of Proposition \ref{propositionsubgaussianoracle}, for any sub-Gaussian vector $\zeta$ with parameter $K_\zeta$ independent of $Z$,
		\begin{enumerate}
			\item $$ \e \left( \sum_{m\in\m_u} w_m \left\Vert \pp{m} \mu \right\Vert^2 \right) = o(n^\tau). $$
			
			\item $$ \e \left( \sum_{m\in\m_u} w_m \mu^\T \pp{m} \zeta \right) = o(n^\tau). $$
		\end{enumerate}
	\end{lemma}
	\begin{proof}
		For $m\in\m_u$, let
		\begin{align*}
			r_m\defined \left\Vert \pp{m} \mu \right\Vert^2.
		\end{align*}
		Fixing $t>0$ arbitrarily, define the set
		\begin{align*}
			\A_t \defined \left\{ m\in\m_u : r_m \leq tn^{\tau} \right\}.
		\end{align*}
		Now,
		\begin{align*}
			\e\left( \sum_{m\in\m_u} w_m r_m \right)
			= \e\left( \sum_{m\in\A_t} w_m r_m \right)
			+ \e\left( \sum_{m\in\A_t^\C} w_m r_m \right)
		\end{align*}
		By the definition of $A_t$,
		\begin{align*}
			\limsup_{n\to\infty} n^{-\tau} \e\left( \sum_{m\in\A_t} w_m r_m \right) \leq t.
		\end{align*}
		For $\A_t^\C$, fix a value of $a>0$, which will be determined later, and define the set $\Tset_a$ as
		\begin{align*}
			\Tset_a \defined \bigcap_{m\in\m_u} \left\{ \left\Vert P_m \xi \right\Vert^2 \leq K_\xi^2 \left( u + 2\sqrt{uan^{\tau}} + 2an^{\tau} \right) \right\}.
		\end{align*}
		By the calculations from equation \eqref{equationtsettcprob}, it follows that
		\begin{align}\label{equationtsetacprob}
			\p \left( \Tset_a^\C \right) \leq \exp\left( -an^\tau + \log\left( |\m_u|\right) \right).
		\end{align}
		Moreover, note that, by assumption,
		\begin{align*}
			\limsup_{n\to\infty} \sup_{m\in\m_u} n^{-1} r_m
			\leq \limsup_{n\to\infty} n^{-1} \left\Vert \mu \right\Vert^2
			\leq C,
		\end{align*}
		for some constant $C>0$.  Then, for $n$ sufficiently large,
		\begin{align}\label{equationatcdecomposition}
			\begin{aligned}
				n^{-\tau} \e \left( \sum_{m\in\A_t^\C} w_m r_m \right)
				\leq 2C n^{1-\tau} \sum_{m\in\A_t^\C} \e \left( w_m \right)
				\leq 2C n^{1-\tau} \sum_{m\in\A_t^\C} \left( \e \left( w_m \indic{\Tset_a} \right) + \p\left( \Tset_a^\C \right) \right).
			\end{aligned}
		\end{align}
		Fix $m\in\A_t^\C$ temporarily and let $S$ be any weakly sparse set for $\mu$.  Then, we have that
		\begin{align*}
			w_m \indic{\Tset_a} &\leq \exp\left( -\frac{1}{\alpha} \left( \left\Vert \pp{m} Y \right\Vert^2 - \left\Vert \pp{S} Y \right\Vert^2 \right) \right)\indic{\Tset_a}\\
			&\leq \exp\left( -\frac{1}{\alpha} \left( r_m - r_S + 2\mu^\T\pp{m}\xi - 2\mu^\T\pp{S}\xi - K_\xi^2 \left( u + 2\sqrt{uan^{\tau}} + 2an^{\tau} \right) \right) \right).
		\end{align*}
		By Cauchy-Schwarz,
		\begin{align*}
			\e\left( w_m \indic{\Tset_a} \right) \leq &\exp\left( -\frac{1}{\alpha}\left( r_m - r_S - K_\xi^2 \left( u + 2\sqrt{uan^{\tau}} + 2an^{\tau} \right) \right) \right)\\
			&\times \left( \e \exp\left( -\frac{4}{\alpha} \mu^\T\pp{m}\xi \right) \right)^{1/2}
			\left( \e \exp\left( \frac{4}{\alpha} \mu^\T\pp{S}\xi \right) \right)^{1/2}.
		\end{align*}
		Computing each of the Laplace transforms directly, it follows that
		\begin{align*}
			\e \exp\left( -\frac{4}{\alpha} \mu^\T\pp{m}\xi \right) \leq \exp\left( \frac{8K_\xi^2}{\alpha^2} r_m \right).
		\end{align*}
		Here, we have used Definition \ref{definitionsg}.  Similarly,
		\begin{align*}
			\e \exp\left( \frac{4}{\alpha} \mu^\T\pp{S}\xi \right) \leq \exp\left( \frac{8K_\xi^2}{\alpha^2} r_S \right).
		\end{align*}
		Hence,
		\begin{align*}
			\e\left( w_m \indic{\Tset_a} \right) &\leq \exp\left( -\frac{1}{\alpha} \left( \left(1 - \frac{4K_\xi^2}{\alpha}\right)r_m - \left( 1 + \frac{4K_\xi^2}{\alpha}\right)r_S - K_\xi^2 \left( u + 2\sqrt{uan^{\tau}} + 2an^{\tau} \right) \right) \right)\\
			&\leq \exp\left( -\frac{1}{\alpha} \left( \left(1 - \frac{4K_\xi^2}{\alpha}\right)tn^\tau - \left( 1 + \frac{4K_\xi^2}{\alpha}\right)r_S - K_\xi^2 \left( u + 2\sqrt{uan^{\tau}} + 2an^{\tau} \right)  \right) \right).
		\end{align*}
		The second inequality follows from the fact that $m\in\A_t^\C$.  Since $u = o(n^\tau/\log(p))$, setting $a < \left(1 - 4K_\xi^2/\alpha\right)t/2$ yields
		\begin{align}\label{equationweightsatcconcentration}
			\e\left( w_m \indic{\Tset_a} \right) \leq \exp\left( -\frac{1}{\alpha}\left(\left( 1 - \frac{4K_\xi^2}{\alpha}\right)t - 2a\right)n^\tau + o(n^\tau) \right)
		\end{align}
		Combining equations \eqref{equationtsetacprob}, \eqref{equationatcdecomposition}, and \eqref{equationweightsatcconcentration}, it follows that
		\begin{align*}
			\limsup_{n\to\infty} n^{-\tau} \e\left( \sum_{m\in\A_t^\C} w_m r_m \right) = 0.
		\end{align*}
		Therefore,
		\begin{align*}
			\limsup_{n\to\infty} n^{-\tau} \e\left( \sum_{m\in\m_u} w_m r_m \right) \leq t.
		\end{align*}
		Since $t>0$ was arbitrary, this proves the first claim.  For the second half, define the set $\Fset_t$ as
		\begin{align*}
			\Fset_t \defined \bigcap_{m\in\A_t} \left\{ \left| \mu^\T \pp{m} \zeta \right| \leq tn^{\tau} \right\}.
		\end{align*}
		For a fixed $m\in\A_t$, it will follow by a Chernoff bound that, for some constant $c>0$,
		\begin{align*}
			\p\left( \left| \mu^\T \pp{m} \zeta \right| > tn^{\tau} \right) \leq 2 \exp \left( -\frac{ct^2 n^{2\tau}}{K_\zeta^2 r_m} \right) \leq 2 \exp \left( -\frac{ct n^{\tau}}{K_\zeta^2} \right).
		\end{align*}
		Therefore, an upper bound for $\p\left( \Fset_t^\C \right)$ is given by
		\begin{align}\label{equationfsettcprob}
			\p\left( \Fset_t^\C \right) \leq 2 \exp\left( -\frac{ct n^{\tau}}{K_\zeta^2} + \log(|\A_t|) \right).
		\end{align}
		Now,
		\begin{align*}
			\e\left(\sum_{m\in\A_t} w_{m} \left|\mu^\T\pp{m}\zeta\right|\right)
			= \e\left(\sum_{m\in\A_t} w_{m} \left|\mu^\T\pp{m}\zeta\right| \indic{\Fset_t} \right)
			+ \e\left(\sum_{m\in\A_t} w_{m} \left|\mu^\T\pp{m}\zeta\right| \indic{\Fset_t^\C} \right).
		\end{align*}
		By the definition of $\Fset_t$, it follows that
		\begin{align*}
			\e\left(\sum_{m\in\A_t} w_{m} \left|\mu^\T\pp{m}\zeta\right| \indic{\Fset_t} \right) \leq tn^\tau.
		\end{align*}
		On $\Fset_t^\C$, two applications of Cauchy-Schwarz and equation \eqref{equationfsettcprob} yields
		\begin{align*}
			\limsup_{n\to\infty} n^{-\tau}&\e\left(\sum_{m\in\A_t} w_{m} \left|\mu^\T\pp{m}\zeta\right| \indic{\Fset_t^\C} \right)\\
			&\leq \limsup_{n\to\infty}n^{-\tau}\left\Vert \mu \right\Vert \e \left( \left\Vert \zeta \right\Vert \indic{\Fset_t^\C} \right)\\
			&\leq \limsup_{n\to\infty}n^{-\tau}\left\Vert \mu \right\Vert \left(\e \left\Vert \zeta \right\Vert^2\right)^{1/2} \left( \p\left( \Fset_t^\C \right) \right)^{1/2}\\
			&= 0.
		\end{align*}
		Furthermore, on $\A_t^\C$, by another two applications of Cauchy-Schwarz,
		\begin{align*}
			\limsup_{n\to\infty} n^{-\tau}&\e\left(\sum_{m\in\A_t^\C} w_{m} \left|\mu^\T\pp{m}\zeta\right|\right)\\
			&\leq \limsup_{n\to\infty} n^{-\tau}\left\Vert \mu \right\Vert \sum_{m\in\A_t^\C} \e\left( w_{m} \left\Vert \zeta \right\Vert\right)\\
			&\leq \limsup_{n\to\infty} n^{-\tau}\left\Vert \mu \right\Vert \sum_{m\in\A_t^\C} \left( \e w_{m}^2 \right)^{1/2} \left(\e \left\Vert \zeta \right\Vert^2\right)^{1/2}\\
			&\leq \limsup_{n\to\infty} n^{-\tau} \left\Vert \mu \right\Vert \left(\e \left\Vert \zeta \right\Vert^2\right)^{1/2} \sum_{m\in\A_t^\C} \left( \e w_{m} \right)^{1/2}\\
			&\leq \limsup_{n\to\infty} n^{-\tau} \left\Vert \mu \right\Vert \left(\e \left\Vert \zeta \right\Vert^2\right)^{1/2} \sum_{m\in\A_t^\C} \left( \e \left( w_{m} \indic{\Tset_a} \right) + \p\left( \Tset_a^\C \right) \right)^{1/2}\\
			&= 0,
		\end{align*}
		where the limit follows by equations \eqref{equationtsetacprob} and \eqref{equationweightsatcconcentration}.  Since $t>0$ was arbitrary, this proves the second claim and finishes the proof.
	\end{proof}
	
	Immediately, we have the following corollary for random designs when the mean vector is assumed to be weakly sparse in probability.
	\begin{corollary}\label{corollarybiasother}
		Consider the setup of Lemma \ref{lemmaq=1nonbias}.  If $\mu$ is weakly sparse relative to $Z$ in probability and $\left\Vert \mu \right\Vert^2 = \Op(n^\tau)$, then
		\begin{enumerate}
			\item $$ \left( \sum_{m\in\m_u} w_m \left\Vert \pp{m} \mu \right\Vert^2 \right) = \op(n^\tau). $$
			
			\item $$ \left( \sum_{m\in\m_u} w_m \mu^\T \pp{m} \zeta \right) = \op(n^\tau). $$
		\end{enumerate}
	\end{corollary}
	
	With these lemmata, we can now prove Proposition \ref{propositionsubgaussianoracle}.
	
	\begin{proof}[Proof of Proposition \ref{propositionsubgaussianoracle}]
		Indeed, by convexity of the norm, it follows that
		\begin{align*}
			\left\Vert \sum_{m\in\m_u} w_m Z\hat{\gamma}_m - \mu \right\Vert^2
			\leq \sum_{m\in\m_u} w_m \left\Vert \pp{m} \mu \right\Vert^2
			+ \sum_{m\in\m_u} w_m \left\Vert P_m \xi \right\Vert^2 .
		\end{align*}
		Applying Lemmata \ref{lemmabiashw} and \ref{lemmabiasother} finishes the proof.
	\end{proof}
	
	Instead of directly proving Theorem \ref{theoremq=1asympdist}, we will decompose the estimator and prove each part separately.
	Indeed, we note that
	\begin{align*}
		\betahatew = \frac{\left( \nu - Z\deltahatew + \eta \right)^\T \left( \mu - Z\thetahatew + \eta\beta + \epsilon \right)}{\left\Vert X - Z\deltahatew \right\Vert^2}.
	\end{align*}
	Then,
	\begin{align*}
		\begin{aligned}
			\sqrt{n}\betahatew = \Big( &\left(\nu - Z\deltahatew \right)^\T \left( \mu - Z\thetahatew + \eta\beta + \epsilon \right) + \eta^\T\left( \mu - Z\thetahatew \right)\\ &+ \eta^\T\eta\beta + \eta^\T\epsilon \Big) \times \frac{1}{\sqrt{n}\sigma_\eta^2} \times \frac{n\sigma_\eta^2}{\left\Vert X - Z\deltahatew \right\Vert^2}.
		\end{aligned}
	\end{align*}
	We will start by proving that the first line, which corresponds to the bias from inexact orthogonalization, converges to zero.
	\begin{lemma}\label{lemmaq=1bias}
		Consider the models given in equations \eqref{equationplmy2} and \eqref{equationplmx}.  Under assumptions \ref{assumptionq=1boundednorms} -- \ref{assumptionq=1sparsity},
		\begin{align*}
			\left(\nu - Z\deltahatew \right)^\T \left( \mu - Z\thetahatew + \eta\beta + \epsilon \right) + \eta^\T\left( \mu - Z\thetahatew \right) = \op(\sqrt{n}).
		\end{align*}
	\end{lemma}
	\begin{proof}
		Without the loss of generality, we will assume that $u \defined u_\theta = u_\delta$.  Expanding, we have
		\begin{align*}
			\left(\nu - Z\deltahatew \right)^\T \left( \mu - Z\thetahatew \right) + \left(\nu - Z\deltahatew \right)^\T\left( \eta\beta + \epsilon \right) + \eta^\T\left( \mu - Z\thetahatew \right).
		\end{align*}
		We will consider each of the three terms separately.  By Cauchy-Schwarz and Corollary \ref{corollarysubgaussianoracle}, it follows that
		\begin{align*}
			\left| \left(\nu - Z\deltahatew \right)^\T \left( \mu - Z\thetahatew \right) \right| \leq \left\Vert \nu - Z\deltahatew \right\Vert \left\Vert \mu - Z\thetahatew \right\Vert = \op(\sqrt{n}).
		\end{align*}
		For the second term, we may further expand to obtain
		\begin{align*}
			\left(\nu - Z\deltahatew \right)^\T\left( \eta\beta + \epsilon \right) =& \sum_{m\in\m_u} w_{m,X} \left( \pp{m} \nu - P_m\eta \right)^\T \left( \eta\beta + \epsilon \right)\\
			=& \sum_{m\in\m_u} w_{m,X} \nu^\T \pp{m} \left( \eta\beta + \epsilon \right) + \frac{1}{2}\sum_{m\in\m_u} w_{m,X} \left\Vert P_{m} \epsilon \right\Vert^2\\
			&- \frac{1}{2}\sum_{m\in\m_u} w_{m,X} \left\Vert P_m \left( \eta + \epsilon \right) \right\Vert^2 \\
			&- \left(\beta - \frac{1}{2}\right) \sum_{m\in\m_u} w_m \left\Vert P_m \eta \right\Vert^2 .
		\end{align*}
		Applying Lemma \ref{lemmabiashw} and Corollary \ref{corollarybiasother}, it follows that
		\begin{align*}
			\left(\nu - Z\deltahatew \right)^\T\left( \eta\beta + \epsilon \right) = \op(\sqrt{n}).
		\end{align*}
		Finally,
		\begin{align*}
			\eta^\T\left( \mu - Z\thetahatew \right) =& \sum_{m\in\m_u} w_{m,Y} \eta^\T \left( \pp{m} \mu - P_m\left( \eta\beta + \epsilon \right) \right)\\
			=& \sum_{m\in\m_u} w_{m,Y} \eta^\T \pp{m} \mu - \frac{1}{2} \sum_{m\in\m_u} w_{m,Y} \left\Vert P_m\left( \eta(\beta+1) + \epsilon \right) \right\Vert^2 \\&+ \frac{1}{2} \sum_{m\in\m_u} w_{m,Y} \left\Vert P_m \left( \eta\beta + \epsilon \right) \right\Vert^2 + \frac{1}{2} \sum_{m\in\m_u} w_{m,Y} \left\Vert P_m \eta \right\Vert^2.
		\end{align*}
		Again, applying Lemma \ref{lemmabiashw} and Corollary \ref{corollarybiasother}, it follows that
		\begin{align*}
			\eta^\T\left( \mu - Z\thetahatew \right) = \op(\sqrt{n}).
		\end{align*}
		This finishes the proof.
	\end{proof}
	\begin{lemma}\label{lemmaq=1nonbias}
		Consider the models given in equations (\ref{equationplmy2}) and (\ref{equationplmx}).  Under assumptions \ref{assumptionq=1boundednorms} -- \ref{assumptionq=1sparsity},
		\begin{enumerate}
			\item $$\sqrt{n}\left( \frac{\eta^\T\eta\beta}{\left\Vert X - Z\deltahatew \right\Vert^2} - \beta \right) \conp 0.$$
			\item $$n^{-1/2}\frac{\eta^\T\epsilon}{\sigma_\eta^2} \cond \n\left( 0 , \frac{\sigma_\epsilon^2}{\sigma_\eta^2}\right).$$
			\item $$\frac{n\sigma_\eta^2}{\left\Vert X - Z\deltahatew \right\Vert^2} \conp 1.$$
		\end{enumerate}
	\end{lemma}
	\begin{proof}
		Indeed, expanding the denominator, we see that
		\begin{align*}
			\left\Vert X - Z\deltahatew \right\Vert^2 = \left\Vert \nu - Z\deltahatew \right\Vert^2 + 2\eta^\T\left( \nu - Z\deltahatew \right) + \left\Vert \eta \right\Vert^2.
		\end{align*}
		By Corollary \ref{corollarysubgaussianoracle} and Lemma \ref{lemmaq=1bias}, it follows that
		\begin{align*}
			\left\Vert X - Z\deltahatew \right\Vert^2 = \op(\sqrt{n}) + \left\Vert \eta \right\Vert^2.
		\end{align*}
		Then, by the Law of Large Numbers,
		\begin{align*}
			\frac{1}{n} \left\Vert X - Z\deltahatew \right\Vert^2 \conp \sigma_\eta^2.
		\end{align*}
		This proves the third claim.  Now, by direct substitution, we have that
		\begin{align*}
			\sqrt{n}\left( \frac{\left(\left\Vert X - Z\deltahatew \right\Vert^2 + \op(\sqrt{n})\right)\beta}{\left\Vert X - Z\deltahatew \right\Vert^2} - \beta \right) = \frac{n}{\left\Vert X - Z\deltahatew \right\Vert^2} \frac{\op(\sqrt{n})}{\sqrt{n}} = \op(1),
		\end{align*}
		which proves the first claim.  The second claim follows by the Central Limit Theorem, which finishes the proof.
	\end{proof}
	
	\begin{proof}[Proof of Theorem \ref{theoremq=1asympdist}]
		The proof follows by combining Lemmata \ref{lemmaq=1bias} and \ref{lemmaq=1nonbias}.
	\end{proof}

	\subsection{Proofs for Section \ref{sectionsparsity}}
	\begin{proof}[Proof of Theorem \ref{theoremjrconjecture}]
		Suppose that $s_\delta = o(\sqrt{n}/\log(p))$.  We will consider a sequence of $\vartheta \in \Theta(s_\gamma,s_\delta)$ such that $S_\gamma \cap S_\delta = \varnothing$ and $\delta \geq 0$ componentwise.  We will construct $\Sigma_{Z,Z}$ implicitly.  For $j\in S_\delta^\C$, let
		\begin{align*}
			Z_j \overset{i.i.d.}{\sim} \n_n \left( 0_n , I_n \right).
		\end{align*}
		Before defining $Z_j$ for $j\in S_\delta$, we will need to define another Gaussian matrix $\Xi \in \R^{n\times p}$.  For $j\in S_\delta^\C$, set $\Xi_j = 0_n$.  Then, for $j\in S_\delta$,
		\begin{align*}
			\Xi_j \overset{i.i.d.}{\sim} \n_n \left( 0_n , \tau_n^2 I_n \right),
		\end{align*}
		independent of $Z_k$ for all $k \in S_\delta^\C$; the value $\tau_n^2 > 0$ will be determined later.  Now, for $j\in S_\delta$, we will let
		\begin{align*}
			Z_j = Z\gamma + \Xi_j.
		\end{align*}
		Therefore, it follows that
		\begin{align*}
			Z\delta = Z\gamma \left\Vert \delta \right\Vert_1 + \Xi\delta.
		\end{align*}
		By a direct calculation,
		\begin{align*}
			\cov( \left(Z\delta\right)_1 , \left(Z\gamma\right)_1) = \cov(\left(Z\gamma\right)_1 \left\Vert \delta \right\Vert_1 + \left(\Xi\delta\right)_1 , \left(Z\gamma\right)_1) = \var(\left(Z\gamma\right)_1)\left\Vert \delta \right\Vert_1.
		\end{align*}
		Moreover,
		\begin{align*}
			\var(\left(Z\delta\right)_1) = \var(\left(Z\gamma\right)_1 \left\Vert \delta \right\Vert_1 + \left(\Xi\delta\right)_1) = \var(\left(Z\gamma\right)_1)\left\Vert \delta \right\Vert_1^2 + \tau_n^2\left\Vert \delta \right\Vert_2^2.
		\end{align*}
		Choosing $\tau_n^2 \to 0$ sufficiently fast, it will follow that
		\begin{align*}
			\var(\left(Z\delta\right)_1) = \var(\left(Z\gamma\right)_1)\left\Vert \delta \right\Vert^2_1 + o\left(n^{-1/2}\right).
		\end{align*}
		Hence, this implies that
		\begin{align*}
			\cov\left( \left(Z\delta\right)_1 , \left( Z\gamma \right)_1 \right) = \sqrt{\var(\left(Z\delta\right)_1)\var(\left(Z\gamma\right)_1)} + o(n^{-1/2}).
		\end{align*}
		Now, note that
		\begin{align*}
			\cov\left( \left(Z\delta\right)_1 , \left( Z\gamma \right)_1 \right) = \cov\left( X_1 , Y_1 \right) - \beta \var\left( X_1 \right).
		\end{align*}
		Let $\betahat$ be any $\sqrt{n}$-consistent estimator for $\beta$.  Then,
		\begin{align*}
			\frac{1}{n} \left( X^\T Y - \betahat X^\T X \right)
		\end{align*}
		is a $\sqrt{n}$-consistent estimator for $\cov\left( \left(Z\delta\right)_1 , \left( Z\gamma \right)_1 \right)$.  Consider an oracle that has access to the set $S_\delta$, knows $S_\delta \cap S_\gamma = \varnothing$, and knows the covariance structure of the design.  Then, since $s_\delta = o(\sqrt{n}/\log(p))$, a $\sqrt{n}$-consistent estimator for $\var(\left(Z\delta\right)_1)$ is given by Theorem \ref{theoremsigmamu}.  This implies that there exists a $\sqrt{n}$-consistent estimator for $\var(\left(Z\gamma\right)_1)$.  By the minimax lower bounds established by \citeA{cai2018}, it follows that, in order to have a $\sqrt{n}$-consistent estimator for $\var(\left(Z\gamma\right)_1)$, it must be the case that $s_\gamma = \mathcal{O}(\sqrt{n}/\log(p))$.  This proves half of the claim.  The other half follows by symmetry, which finishes the proof.
	\end{proof}

	\bibliographystyle{newapa}
	\bibliography{NoDataSplittingref}

	\newpage
	

	\title{Supplement to ``Inference Without Compatibility''}
	
	\author{Michael Law \and Ya\hspace{-.1em}'\hspace{-.1em}acov Ritov}
	\date{
		University of Michigan\\
		\today
	}
	
	\maketitle
	
	\begin{abstract}
		This supplement provides additional simulation results as well as the proofs for Sections 2.2, 2.3, 3.1, and 3.2.
	\end{abstract}
	
	\newcounter{suppsection}
	\addtocounter{suppsection}{1}
	\def\thesection{S\arabic{suppsection}}
	
	\section{Additional Simulation Results}
	In this section, we include additional results for the simulations of Section \ref{sectionsimulations}.
	
	\def\thetable{S\arabic{table}}
	
\begin{table}[]
\centering
\caption{Simulations for $\beta$ with Gaussian design and errors when q=3 and $\beta=$0} 
\label{tableq3beta0desz}
\begin{tabular}{|l|l|rrrrrrrr|}
   \hline
 & $snr_X$ & 2 & 2 & 2 & 2 & 1000 & 1000 & 1000 & 1000 \\ 
   & $\rho$ & 0 & 0 & 0.8 & 0.8 & 0 & 0 & 0.8 & 0.8 \\ 
   & $s_\delta,s_\gamma$ & 3 & 15 & 3 & 15 & 3 & 15 & 3 & 15 \\ 
   \hline
 & LS & 0.924 & 0.886 & 0.922 & 0.936 & 0.910 & 0.872 & 0.942 & 0.946 \\ 
   & SILM & 0.936 & 0.894 & 0.960 & 0.964 & 0.940 & 0.788 & 0.896 & 0.874 \\ 
  AveCov & $\text{EW}_{I}$ & 0.944 & 0.886 & 0.956 & 0.950 & 0.952 & 0.800 & 0.976 & 0.978 \\ 
   & $\text{EW}_{II}$ & 0.978 & 0.942 & 0.978 & 0.976 & 0.976 & 0.908 & 0.992 & 0.994 \\ 
   & $\text{EW}_{III}$ & 0.992 & 0.964 & 0.988 & 0.990 & 0.990 & 0.952 & 0.998 & 0.994 \\ 
   \hline
\end{tabular}
\end{table}

\begin{table}[]
\centering
\caption{Simulations for $\beta$ with Gaussian design and errors when q=1 and $\beta=$1} 
\label{tableq1beta1desz}
\begin{tabular}{|l|l|rrrrrrrr|}
   \hline
 & $snr_X$ & 2 & 2 & 2 & 2 & 1000 & 1000 & 1000 & 1000 \\ 
   & $\rho$ & 0 & 0 & 0.8 & 0.8 & 0 & 0 & 0.8 & 0.8 \\ 
   & $s_\delta,s_\gamma$ & 3 & 15 & 3 & 15 & 3 & 15 & 3 & 15 \\ 
   \hline
 & LS & 0.946 & 0.878 & 0.938 & 0.946 & 0.934 & 0.900 & 0.948 & 0.944 \\ 
   & DLA & 0.928 & 0.880 & 0.934 & 0.946 & 0.904 & 0.856 & 0.352 & 0.238 \\ 
   & SILM & 0.932 & 0.872 & 0.936 & 0.956 & 0.918 & 0.858 & 0.130 & 0.034 \\ 
  AvgCov & DML & 1.000 & 0.996 & 0.998 & 0.994 & 0.990 & 0.982 & 1.000 & 1.000 \\ 
   & $\text{EW}_{I}$ & 0.830 & 0.768 & 0.922 & 0.932 & 0.932 & 0.862 & 0.976 & 0.984 \\ 
   & $\text{EW}_{II}$ & 0.866 & 0.810 & 0.952 & 0.962 & 0.956 & 0.900 & 0.984 & 0.992 \\ 
   & $\text{EW}_{III}$ & 0.904 & 0.852 & 0.974 & 0.976 & 0.964 & 0.922 & 0.990 & 0.998 \\ 
   \hline
 & LS & 0.428 & 0.463 & 0.591 & 0.688 & 0.429 & 0.466 & 0.932 & 1.450 \\ 
   & DLA & 0.502 & 0.539 & 0.693 & 0.699 & 0.539 & 0.555 & 0.548 & 0.506 \\ 
   & SILM & 0.549 & 0.579 & 0.683 & 0.709 & 0.647 & 0.636 & 0.673 & 0.640 \\ 
  AvgLen & DML & 1.190 & 1.180 & 1.180 & 1.160 & 2.800 & 1.670 & 17.300 & 17.200 \\ 
   & $\text{EW}_{I}$ & 0.640 & 0.655 & 0.711 & 0.717 & 1.080 & 0.805 & 1.910 & 1.870 \\ 
   & $\text{EW}_{II}$ & 0.696 & 0.717 & 0.778 & 0.801 & 1.180 & 0.883 & 2.130 & 2.120 \\ 
   & $\text{EW}_{III}$ & 0.746 & 0.773 & 0.839 & 0.877 & 1.270 & 0.953 & 2.320 & 2.350 \\ 
   \hline
\end{tabular}
\end{table}

\begin{table}[]
\centering
\caption{Simulations for $\beta$ with Gaussian design and errors when q=3 and $\beta=$1} 
\label{tableq3beta1desz}
\begin{tabular}{|l|l|rrrrrrrr|}
   \hline
 & $snr_X$ & 2 & 2 & 2 & 2 & 1000 & 1000 & 1000 & 1000 \\ 
   & $\rho$ & 0 & 0 & 0.8 & 0.8 & 0 & 0 & 0.8 & 0.8 \\ 
   & $s_\delta,s_\gamma$ & 3 & 15 & 3 & 15 & 3 & 15 & 3 & 15 \\ 
   \hline
 & LS & 0.914 & 0.882 & 0.928 & 0.938 & 0.950 & 0.848 & 0.924 & 0.948 \\ 
   & SILM & 0.862 & 0.684 & 0.886 & 0.916 & 0.892 & 0.644 & 0.020 & 0.004 \\ 
  AveCov & $\text{EW}_{I}$ & 0.522 & 0.636 & 0.766 & 0.838 & 0.832 & 0.736 & 0.858 & 0.928 \\ 
   & $\text{EW}_{II}$ & 0.568 & 0.696 & 0.828 & 0.886 & 0.876 & 0.790 & 0.908 & 0.968 \\ 
   & $\text{EW}_{III}$ & 0.630 & 0.746 & 0.858 & 0.938 & 0.900 & 0.828 & 0.952 & 0.978 \\ 
   \hline
\end{tabular}
\end{table}

\begin{table}[]
\centering
\caption{Simulations for $\beta$ with double exponential design and errors when q=1 and $\beta=$0} 
\label{tableq1beta0dese}
\begin{tabular}{|l|l|rrrrrrrr|}
   \hline
 & $snr_X$ & 2 & 2 & 2 & 2 & 1000 & 1000 & 1000 & 1000 \\ 
   & $\rho$ & 0 & 0 & 0.8 & 0.8 & 0 & 0 & 0.8 & 0.8 \\ 
   & $s_\delta,s_\gamma$ & 3 & 15 & 3 & 15 & 3 & 15 & 3 & 15 \\ 
   \hline
 & LS & 0.942 & 0.900 & 0.966 & 0.954 & 0.950 & 0.892 & 0.940 & 0.946 \\ 
   & DLA & 0.960 & 0.876 & 0.974 & 0.964 & 0.950 & 0.872 & 0.154 & 0.102 \\ 
   & SILM & 0.954 & 0.876 & 0.968 & 0.960 & 0.954 & 0.838 & 0.892 & 0.868 \\ 
  AvgCov & DML & 0.960 & 0.878 & 0.930 & 0.914 & 0.980 & 0.848 & 1.000 & 1.000 \\ 
   & $\text{EW}_{I}$ & 0.950 & 0.858 & 0.962 & 0.962 & 0.956 & 0.866 & 0.960 & 0.954 \\ 
   & $\text{EW}_{II}$ & 0.972 & 0.918 & 0.978 & 0.976 & 0.976 & 0.914 & 0.970 & 0.976 \\ 
   & $\text{EW}_{III}$ & 0.980 & 0.938 & 0.990 & 0.984 & 0.990 & 0.950 & 0.980 & 0.990 \\ 
   \hline
 & LS & 0.456 & 0.492 & 0.683 & 0.829 & 0.432 & 0.466 & 0.910 & 1.450 \\ 
   & DLA & 0.534 & 0.565 & 0.813 & 0.821 & 0.530 & 0.545 & 0.520 & 0.490 \\ 
   & SILM & 0.574 & 0.596 & 0.785 & 0.825 & 0.619 & 0.603 & 0.611 & 0.592 \\ 
  AvgLen & DML & 0.756 & 0.702 & 0.874 & 0.875 & 1.480 & 0.892 & 12.700 & 13.300 \\ 
   & $\text{EW}_{I}$ & 0.716 & 0.691 & 0.856 & 0.877 & 1.070 & 0.798 & 2.000 & 1.910 \\ 
   & $\text{EW}_{II}$ & 0.792 & 0.772 & 0.932 & 0.974 & 1.180 & 0.891 & 2.180 & 2.120 \\ 
   & $\text{EW}_{III}$ & 0.860 & 0.844 & 1.000 & 1.060 & 1.280 & 0.973 & 2.330 & 2.300 \\ 
   \hline
\end{tabular}
\end{table}

\begin{table}[]
\centering
\caption{Simulations for $\beta$ with double exponential design and errors when q=3 and $\beta=$0} 
\label{tableq3beta0dese}
\begin{tabular}{|l|l|rrrrrrrr|}
   \hline
 & $snr_X$ & 2 & 2 & 2 & 2 & 1000 & 1000 & 1000 & 1000 \\ 
   & $\rho$ & 0 & 0 & 0.8 & 0.8 & 0 & 0 & 0.8 & 0.8 \\ 
   & $s_\delta,s_\gamma$ & 3 & 15 & 3 & 15 & 3 & 15 & 3 & 15 \\ 
   \hline
 & LS & 0.928 & 0.834 & 0.904 & 0.930 & 0.940 & 0.884 & 0.926 & 0.940 \\ 
   & SILM & 0.954 & 0.856 & 0.956 & 0.966 & 0.950 & 0.782 & 0.874 & 0.858 \\ 
  AveCov & $\text{EW}_{I}$ & 0.958 & 0.850 & 0.954 & 0.958 & 0.946 & 0.794 & 0.966 & 0.976 \\ 
   & $\text{EW}_{II}$ & 0.984 & 0.936 & 0.990 & 0.978 & 0.972 & 0.886 & 0.990 & 0.992 \\ 
   & $\text{EW}_{III}$ & 0.992 & 0.966 & 0.992 & 0.992 & 0.984 & 0.932 & 0.990 & 0.996 \\ 
   \hline
\end{tabular}
\end{table}

\begin{table}[]
\centering
\caption{Simulations for $\beta$ with double exponential design and errors when q=1 and $\beta=$1} 
\label{tableq1beta1dese}
\begin{tabular}{|l|l|rrrrrrrr|}
   \hline
 & $snr_X$ & 2 & 2 & 2 & 2 & 1000 & 1000 & 1000 & 1000 \\ 
   & $\rho$ & 0 & 0 & 0.8 & 0.8 & 0 & 0 & 0.8 & 0.8 \\ 
   & $s_\delta,s_\gamma$ & 3 & 15 & 3 & 15 & 3 & 15 & 3 & 15 \\ 
   \hline
 & LS & 0.934 & 0.902 & 0.946 & 0.930 & 0.940 & 0.914 & 0.934 & 0.954 \\ 
   & DLA & 0.940 & 0.884 & 0.920 & 0.918 & 0.910 & 0.876 & 0.332 & 0.280 \\ 
   & SILM & 0.946 & 0.892 & 0.910 & 0.904 & 0.922 & 0.842 & 0.082 & 0.014 \\ 
  AvgCov & DML & 0.998 & 1.000 & 0.994 & 0.994 & 0.986 & 0.956 & 1.000 & 1.000 \\ 
   & $\text{EW}_{I}$ & 0.894 & 0.818 & 0.910 & 0.942 & 0.934 & 0.860 & 0.958 & 0.974 \\ 
   & $\text{EW}_{II}$ & 0.920 & 0.856 & 0.938 & 0.968 & 0.958 & 0.880 & 0.972 & 0.986 \\ 
   & $\text{EW}_{III}$ & 0.942 & 0.886 & 0.956 & 0.974 & 0.974 & 0.902 & 0.988 & 0.990 \\ 
   \hline
 & LS & 0.454 & 0.495 & 0.677 & 0.831 & 0.434 & 0.468 & 0.899 & 1.470 \\ 
   & DLA & 0.541 & 0.573 & 0.803 & 0.829 & 0.545 & 0.550 & 0.519 & 0.491 \\ 
   & SILM & 0.602 & 0.629 & 0.788 & 0.842 & 0.652 & 0.627 & 0.603 & 0.586 \\ 
  AvgLen & DML & 1.360 & 1.270 & 1.480 & 1.460 & 2.830 & 1.660 & 22.200 & 22.000 \\ 
   & $\text{EW}_{I}$ & 0.726 & 0.715 & 0.854 & 0.892 & 1.100 & 0.806 & 2.020 & 1.910 \\ 
   & $\text{EW}_{II}$ & 0.791 & 0.789 & 0.941 & 1.000 & 1.190 & 0.884 & 2.250 & 2.170 \\ 
   & $\text{EW}_{III}$ & 0.851 & 0.855 & 1.020 & 1.100 & 1.280 & 0.953 & 2.460 & 2.390 \\ 
   \hline
\end{tabular}
\end{table}

\begin{table}[]
\centering
\caption{Simulations for $\beta$ with double exponential design and errors when q=3 and $\beta=$1} 
\label{tableq3beta1dese}
\begin{tabular}{|l|l|rrrrrrrr|}
   \hline
 & $snr_X$ & 2 & 2 & 2 & 2 & 1000 & 1000 & 1000 & 1000 \\ 
   & $\rho$ & 0 & 0 & 0.8 & 0.8 & 0 & 0 & 0.8 & 0.8 \\ 
   & $s_\delta,s_\gamma$ & 3 & 15 & 3 & 15 & 3 & 15 & 3 & 15 \\ 
   \hline
 & LS & 0.936 & 0.878 & 0.950 & 0.940 & 0.948 & 0.854 & 0.934 & 0.942 \\ 
   & SILM & 0.864 & 0.682 & 0.876 & 0.906 & 0.878 & 0.624 & 0.010 & 0.002 \\ 
  AveCov & $\text{EW}_{I}$ & 0.582 & 0.692 & 0.818 & 0.844 & 0.828 & 0.714 & 0.868 & 0.950 \\ 
   & $\text{EW}_{II}$ & 0.658 & 0.760 & 0.870 & 0.904 & 0.880 & 0.774 & 0.920 & 0.976 \\ 
   & $\text{EW}_{III}$ & 0.714 & 0.822 & 0.914 & 0.936 & 0.912 & 0.820 & 0.960 & 0.986 \\ 
   \hline
\end{tabular}
\end{table}

\begin{table}[]
\centering
\caption{Simulations for $\beta$ with scaled t design and errors when q=1 and $\beta=$0} 
\label{tableq1beta0dest}
\begin{tabular}{|l|l|rrrrrrrr|}
   \hline
 & $snr_X$ & 2 & 2 & 2 & 2 & 1000 & 1000 & 1000 & 1000 \\ 
   & $\rho$ & 0 & 0 & 0.8 & 0.8 & 0 & 0 & 0.8 & 0.8 \\ 
   & $s_\delta,s_\gamma$ & 3 & 15 & 3 & 15 & 3 & 15 & 3 & 15 \\ 
   \hline
 & LS & 0.958 & 0.910 & 0.942 & 0.938 & 0.956 & 0.904 & 0.950 & 0.960 \\ 
   & DLA & 0.954 & 0.878 & 0.962 & 0.948 & 0.946 & 0.878 & 0.198 & 0.114 \\ 
   & SILM & 0.968 & 0.882 & 0.968 & 0.960 & 0.946 & 0.838 & 0.866 & 0.834 \\ 
  AvgCov & DML & 0.980 & 0.868 & 0.920 & 0.880 & 0.976 & 0.822 & 0.998 & 0.998 \\ 
   & $\text{EW}_{I}$ & 0.950 & 0.846 & 0.956 & 0.956 & 0.966 & 0.816 & 0.968 & 0.974 \\ 
   & $\text{EW}_{II}$ & 0.972 & 0.902 & 0.984 & 0.974 & 0.976 & 0.880 & 0.986 & 0.982 \\ 
   & $\text{EW}_{III}$ & 0.988 & 0.940 & 0.984 & 0.980 & 0.982 & 0.904 & 0.990 & 0.994 \\ 
   \hline
 & LS & 0.490 & 0.515 & 0.750 & 0.933 & 0.453 & 0.479 & 0.951 & 1.600 \\ 
   & DLA & 0.559 & 0.591 & 0.886 & 0.892 & 0.525 & 0.547 & 0.569 & 0.544 \\ 
   & SILM & 0.611 & 0.620 & 0.884 & 0.928 & 0.618 & 0.607 & 0.687 & 0.672 \\ 
  AvgLen & DML & 0.819 & 0.751 & 1.140 & 1.170 & 1.410 & 0.902 & 13.200 & 14.400 \\ 
   & $\text{EW}_{I}$ & 0.806 & 0.739 & 0.967 & 0.982 & 1.100 & 0.817 & 2.240 & 2.110 \\ 
   & $\text{EW}_{II}$ & 0.882 & 0.828 & 1.060 & 1.090 & 1.210 & 0.914 & 2.450 & 2.330 \\ 
   & $\text{EW}_{III}$ & 0.952 & 0.907 & 1.140 & 1.180 & 1.300 & 0.999 & 2.640 & 2.530 \\ 
   \hline
\end{tabular}
\end{table}

\begin{table}[]
\centering
\caption{Simulations for $\beta$ with scaled t design and errors when q=3 and $\beta=$0} 
\label{tableq3beta0dest}
\begin{tabular}{|l|l|rrrrrrrr|}
   \hline
 & $snr_X$ & 2 & 2 & 2 & 2 & 1000 & 1000 & 1000 & 1000 \\ 
   & $\rho$ & 0 & 0 & 0.8 & 0.8 & 0 & 0 & 0.8 & 0.8 \\ 
   & $s_\delta,s_\gamma$ & 3 & 15 & 3 & 15 & 3 & 15 & 3 & 15 \\ 
   \hline
 & LS & 0.922 & 0.872 & 0.940 & 0.944 & 0.936 & 0.856 & 0.926 & 0.936 \\ 
   & SILM & 0.954 & 0.866 & 0.964 & 0.972 & 0.954 & 0.822 & 0.846 & 0.796 \\ 
  AveCov & $\text{EW}_{I}$ & 0.958 & 0.832 & 0.958 & 0.968 & 0.950 & 0.798 & 0.962 & 0.972 \\ 
   & $\text{EW}_{II}$ & 0.980 & 0.924 & 0.978 & 0.990 & 0.986 & 0.886 & 0.986 & 0.994 \\ 
   & $\text{EW}_{III}$ & 0.990 & 0.958 & 0.988 & 0.994 & 0.990 & 0.916 & 0.996 & 0.996 \\ 
   \hline
\end{tabular}
\end{table}

\begin{table}[]
\centering
\caption{Simulations for $\beta$ with scaled t design and errors when q=1 and $\beta=$1} 
\label{tableq1beta1dest}
\begin{tabular}{|l|l|rrrrrrrr|}
   \hline
 & $snr_X$ & 2 & 2 & 2 & 2 & 1000 & 1000 & 1000 & 1000 \\ 
   & $\rho$ & 0 & 0 & 0.8 & 0.8 & 0 & 0 & 0.8 & 0.8 \\ 
   & $s_\delta,s_\gamma$ & 3 & 15 & 3 & 15 & 3 & 15 & 3 & 15 \\ 
   \hline
 & LS & 0.938 & 0.894 & 0.926 & 0.948 & 0.946 & 0.894 & 0.952 & 0.952 \\ 
   & DLA & 0.924 & 0.888 & 0.908 & 0.936 & 0.900 & 0.874 & 0.412 & 0.348 \\ 
   & SILM & 0.886 & 0.934 & 0.900 & 0.914 & 0.908 & 0.842 & 0.112 & 0.046 \\ 
  AvgCov & DML & 0.998 & 1.000 & 0.992 & 0.998 & 0.984 & 0.978 & 1.000 & 1.000 \\ 
   & $\text{EW}_{I}$ & 0.882 & 0.790 & 0.940 & 0.954 & 0.926 & 0.816 & 0.976 & 0.978 \\ 
   & $\text{EW}_{II}$ & 0.920 & 0.836 & 0.964 & 0.978 & 0.956 & 0.868 & 0.982 & 0.990 \\ 
   & $\text{EW}_{III}$ & 0.942 & 0.868 & 0.972 & 0.990 & 0.968 & 0.890 & 0.988 & 0.994 \\ 
   \hline
 & LS & 0.487 & 0.521 & 0.742 & 0.938 & 0.458 & 0.478 & 0.956 & 1.610 \\ 
   & DLA & 0.565 & 0.608 & 0.892 & 0.910 & 0.549 & 0.558 & 0.573 & 0.548 \\ 
   & SILM & 0.643 & 0.657 & 0.902 & 0.957 & 0.659 & 0.638 & 0.696 & 0.668 \\ 
  AvgLen & DML & 1.540 & 1.370 & 1.850 & 1.750 & 2.870 & 1.730 & 22.800 & 23.200 \\ 
   & $\text{EW}_{I}$ & 0.814 & 0.753 & 0.983 & 1.010 & 1.110 & 0.831 & 2.140 & 2.150 \\ 
   & $\text{EW}_{II}$ & 0.886 & 0.830 & 1.090 & 1.140 & 1.210 & 0.911 & 2.390 & 2.420 \\ 
   & $\text{EW}_{III}$ & 0.952 & 0.898 & 1.190 & 1.250 & 1.300 & 0.983 & 2.610 & 2.670 \\ 
   \hline
\end{tabular}
\end{table}
	
\begin{table}[]
\centering
\caption{Simulations for $\beta$ with scaled t design and errors when q=3 and $\beta=$1} 
\label{tableq3beta1dest}
\begin{tabular}{|l|l|rrrrrrrr|}
   \hline
 & $snr_X$ & 2 & 2 & 2 & 2 & 1000 & 1000 & 1000 & 1000 \\ 
   & $\rho$ & 0 & 0 & 0.8 & 0.8 & 0 & 0 & 0.8 & 0.8 \\ 
   & $s_\delta,s_\gamma$ & 3 & 15 & 3 & 15 & 3 & 15 & 3 & 15 \\ 
   \hline
 & LS & 0.950 & 0.906 & 0.944 & 0.918 & 0.918 & 0.874 & 0.950 & 0.948 \\ 
   & SILM & 0.894 & 0.704 & 0.890 & 0.874 & 0.874 & 0.746 & 0.052 & 0.008 \\ 
  AveCov & $\text{EW}_{I}$ & 0.638 & 0.644 & 0.820 & 0.832 & 0.840 & 0.694 & 0.926 & 0.940 \\ 
   & $\text{EW}_{II}$ & 0.682 & 0.716 & 0.904 & 0.906 & 0.890 & 0.746 & 0.962 & 0.984 \\ 
   & $\text{EW}_{III}$ & 0.724 & 0.766 & 0.944 & 0.952 & 0.908 & 0.792 & 0.974 & 0.994 \\ 
   \hline
\end{tabular}
\end{table}

\begin{table}[H]
\centering
\caption{Simulations for $\sigma_{\mu}^2$ with $s_{\gamma}=$15} 
\label{tablemusgamma15}
\begin{tabular}{|l|l|rrrrrr|}
   \hline
 & Distribution & z & z & e & e & t & t \\ 
   & $\rho$ & 0 & 0.8 & 0 & 0.8 & 0 & 0.8 \\ 
   \hline
 & LS & 0.762 & 0.734 & 0.768 & 0.816 & 0.892 & 0.906 \\ 
   & $\text{CHIVE}_{0}$ & 0.134 & 0.492 & 0.152 & 0.464 & 0.228 & 0.460 \\ 
   & $\text{CHIVE}_{2}$ & 0.380 & 0.584 & 0.392 & 0.560 & 0.408 & 0.554 \\ 
  AvgCov & $\text{CHIVE}_{4}$ & 0.514 & 0.676 & 0.554 & 0.638 & 0.540 & 0.674 \\ 
   & $\text{CHIVE}_{6}$ & 0.646 & 0.740 & 0.632 & 0.698 & 0.624 & 0.690 \\ 
   & $\text{EW}_{I}$ & 0.328 & 0.696 & 0.422 & 0.732 & 0.388 & 0.652 \\ 
   & $\text{EW}_{II}$ & 0.628 & 0.756 & 0.630 & 0.784 & 0.588 & 0.786 \\ 
   & $\text{EW}_{III}$ & 0.690 & 0.606 & 0.672 & 0.690 & 0.718 & 0.816 \\ 
   \hline
 & LS & 1.450 & 1.440 & 1.500 & 1.850 & 1.990 & 3.080 \\ 
   & $\text{CHIVE}_{0}$ & 0.538 & 0.873 & 0.583 & 0.999 & 1.060 & 2.140 \\ 
   & $\text{CHIVE}_{2}$ & 1.410 & 1.550 & 1.420 & 1.670 & 1.810 & 2.700 \\ 
  AvgLen & $\text{CHIVE}_{4}$ & 1.910 & 1.980 & 1.910 & 2.090 & 2.270 & 3.070 \\ 
   & $\text{CHIVE}_{6}$ & 2.310 & 2.320 & 2.300 & 2.440 & 2.630 & 3.380 \\ 
   & $\text{EW}_{I}$ & 1.200 & 1.370 & 1.260 & 1.640 & 1.720 & 2.830 \\ 
   & $\text{EW}_{II}$ & 1.280 & 1.340 & 1.320 & 1.630 & 1.790 & 2.820 \\ 
   & $\text{EW}_{III}$ & 1.230 & 1.250 & 1.270 & 1.560 & 1.770 & 2.780 \\ 
   \hline
\end{tabular}
\end{table}

\begin{table}[H]
\centering
\caption{Simulations for $\sigma_{\epsilon}^2$ with $s_{\gamma}=$15} 
\label{tableepsilonsgamma15}
\begin{tabular}{|l|l|rrrrrr|}
   \hline
 & Distribution & z & z & e & e & t & t \\ 
   & $\rho$ & 0 & 0.8 & 0 & 0.8 & 0 & 0.8 \\ 
   \hline
 & LS & 0.874 & 0.864 & 0.870 & 0.848 & 0.876 & 0.864 \\ 
   & SL & 0.308 & 0.646 & 0.386 & 0.620 & 0.466 & 0.606 \\ 
  AvgCov & RCV-SIS & 0.004 & 0.238 & 0.006 & 0.256 & 0.012 & 0.254 \\ 
   & $\text{EW}_{I}$ & 0.514 & 0.630 & 0.554 & 0.650 & 0.532 & 0.648 \\ 
   & $\text{EW}_{II}$ & 0.026 & 0.362 & 0.042 & 0.358 & 0.058 & 0.376 \\ 
   & $\text{EW}_{III}$ & 0.000 & 0.110 & 0.006 & 0.092 & 0.002 & 0.126 \\ 
   \hline
 & LS & 0.481 & 0.462 & 0.467 & 0.478 & 0.479 & 0.483 \\ 
   & SL & 0.781 & 0.702 & 0.766 & 0.722 & 0.753 & 0.717 \\ 
  AvgLen & RCV-SIS & 1.030 & 0.711 & 1.030 & 0.746 & 1.210 & 0.724 \\ 
   & $\text{EW}_{I}$ & 0.613 & 0.498 & 0.600 & 0.515 & 0.589 & 0.507 \\ 
   & $\text{EW}_{II}$ & 0.676 & 0.536 & 0.660 & 0.553 & 0.644 & 0.541 \\ 
   & $\text{EW}_{III}$ & 0.793 & 0.605 & 0.770 & 0.623 & 0.746 & 0.604 \\ 
   \hline
\end{tabular}
\end{table}

	\addtocounter{suppsection}{1}
	
	\section{Proofs}

	\subsection{Proofs for Section \ref{sectionq=1correlated}}
	
	\begin{lemmaS}\label{lemmaq=1correlated}
		Consider the models given in equations (\ref{equationplmy2}) and (\ref{equationplmx}).  Under assumptions \ref{assumptionq=1boundednorms}, \ref{assumptionq=1errordistcorrelated}, and \ref{assumptionq=1sparsity},
		\begin{align*}
			\frac{\eta^\T \epsilon}{\sigma_\eta \Tr(\Sigma_\epsilon)} \cond \n\left( 0 , 1 \right).
		\end{align*}
	\end{lemmaS}
	\begin{proof}
		By the Spectral Theorem, there exists a unitary matrix $\Gamma$ and a diagonal matrix $D$ such that $\Sigma_\epsilon = \Gamma D \Gamma^\T$.  Since $\epsilon$ and $\eta$ are both Gaussian and independent, there exists Gaussian vectors $\zeta\sim\n_n\left( 0 , I_n \right)$ and $\xi \sim \n_n\left( 0_n , I_n \right)$ such that
		\begin{align*}
			\eta^\T \epsilon \equid \sigma_\eta \zeta^\T D^{1/2} \xi.
		\end{align*}
		Then, by the Lindeberg Central Limit Theorem, it follows that
		\begin{align*}
			\frac{\zeta^\T D^{1/2} \xi}{\sqrt{\Tr(D)}} \cond \n\left( 0 , 1 \right).
		\end{align*}
		Noting that $\Tr(D) = \Tr(\Sigma)$ finishes the proof.
	\end{proof}
	\begin{proof}[Proof of Theorem \ref{theoremq=1asympdistcorrelated}]
		The proof follows by combining Lemmata \ref{lemmaq=1bias}, \ref{lemmaq=1nonbias}, and \ref{lemmaq=1correlated}.
	\end{proof}

		\subsection{Proofs for Section \ref{sectionq>1}}
	Similar to the setting where $q=1$, we will proceed in a few stages.
	
	\begin{lemmaS}\label{lemmaq>1nonbias}
		Consider the models given in equations (\ref{equationplmy2}) and (\ref{equationplmx}).  Under assumptions \ref{assumptionq>1boundednorms} -- \ref{assumptionq>1sparsity},
		\begin{enumerate}
			\item $$\left\Vert \left( N - Z\Deltahatew \right)^\T \left( N - Z\Deltahatew \right) \right\Vert = \op(\sqrt{n}).$$
			\item $$\left\Vert \left( N - Z\Deltahatew \right)^\T H \right\Vert = \op(\sqrt{n}).$$
			\item $$n\left(\left( X - Z\Deltahatew \right)^\T \left( X - Z\Deltahatew \right)\right)^{-1} \conp \Sigma_H^{-1}.$$
		\end{enumerate}
	\end{lemmaS}
	\begin{proof}
		Indeed, note that $\left( N - Z\Deltahatew \right)^\T \left( N - Z\Deltahatew \right)$ is a positive definite matrix.  By $q$ applications of Lemma \ref{lemmaq=1bias}, each diagonal element is $\op(\sqrt{n})$, which proves the first claim.  For the second part, Lemma \ref{lemmaq=1bias} again shows that each diagonal element is $\op(\sqrt{n})$.  It is left to show that each off diagonal element is also $\op(\sqrt{n})$.  By symmetry, it suffices to consider the $(1,2)$ element of $\left( N - Z\Deltahatew \right)^\T H$.  For simplicity, we will write $\nu$ to denote the first column of $N$, $\deltahatew$ to denote the first column of $\Deltahatew$, $w_m$ to denote the exponential weights of $\deltahatew$, $\eta$ to denote the first column of $H$, and $\xi$ to denote the second column of $H$.  Then, the $(1,2)$ element can be expressed as
		\begin{align*}
			\left( \nu - Z\deltahatew \right)^\T \xi =& \sum_{m\in\m_u} w_m \nu^\T \pp{m} \xi - \sum_{m\in\m_u} w_m \eta^\T P_m \xi \\
			=&  \sum_{m\in\m_u} w_m \nu^\T \pp{m} \xi - \frac{1}{2} \sum_{m\in\m_u} w_m \left\Vert P_m \left( \xi + \eta \right) \right\Vert^2\\
			&+ \frac{1}{2} \sum_{m\in\m_u} w_m \left\Vert P_m \xi \right\Vert^2 + \frac{1}{2} \sum_{m\in\m_u} w_m \left\Vert P_m \eta \right\Vert^2.
		\end{align*}
		Applying Lemma \ref{lemmabiashw} and Corollary \ref{corollarybiasother} proves the second claim.  Finally, note that
		\begin{align*}
			\Bigg\Vert \Big( X& - Z\Deltahatew \Big)^\T \left( X - Z\Deltahatew \right) - n\Sigma_H \Bigg\Vert\\
			\leq& \left\Vert \left( N - Z\Deltahatew \right)^\T \left( N - Z\Deltahatew \right) \right\Vert +
			2\left\Vert \left( N - Z\Deltahatew \right)^\T H \right\Vert\\& + \left\Vert H^\T H - n\Sigma_H \right\Vert.
		\end{align*}
		We have already shown that the first two terms are $\op(\sqrt{n})$.  For the last term, by the Law of Large Numbers, it follows that
		\begin{align*}
			\left\Vert H^\T H - n\Sigma_H \right\Vert = \op(n).
		\end{align*}
		Therefore,
		\begin{align*}
			\frac{1}{n} \left( X - Z\Deltahatew \right)^\T \left( X - Z\Deltahatew \right) \conp \Sigma_H.
		\end{align*}
		Since $\Sigma_H$ is assumed to be invertible, applying the Continuous Mapping Theorem finishes the proof.
	\end{proof}
	\begin{proof}[Proof of Theorem \ref{theoremq>1asympdist}]
		For convenience, define the following matrices
		\begin{align*}
			&A \defined \left(\left( X - Z\Deltahatew \right)^\T \left( X - Z\Deltahatew \right)\right),\\
			&B \defined \left( N - Z\Deltahatew \right)^\T H,\\
			&C \defined \left( N - Z\Deltahatew \right)^\T \left( N - Z\Deltahatew \right).
		\end{align*}
		Applying Lemma \ref{lemmaq=1bias} to each row separately, we see that
		\begin{align*}
			\betahatew &= \sqrt{n} A^{-1}\left( X - Z\Deltahatew \right)^\T \left( Y - Z\thetahatew \right)\\ &= \sqrt{n} A^{-1} \left(H^\T H\beta + H^\T \epsilon + R\right),
		\end{align*}
		where $\left\Vert R \right\Vert_1 = \op(\sqrt{n})$.  But, from Lemma \ref{lemmaq>1nonbias}, we have that
		\begin{align*}
			n \left\Vert A^{-1} \right\Vert \conp \left\Vert \Sigma_H^{-1} \right\Vert,
		\end{align*}
		which is finite since $\Sigma_H$ is invertible by assumption.  Therefore,
		\begin{align*}
			\left\Vert \sqrt{n} A^{-1} R \right\Vert \leq \left(n \left\Vert A^{-1} \right\Vert\right) \left( n^{-1/2} \left\Vert R \right\Vert_1 \right) \conp 0.
		\end{align*}
		Now, note that
		\begin{align*}
			H^\T H =& A - B - B^\T - C.
		\end{align*}
		Hence,
		\begin{align*}
			\sqrt{n}A^{-1} H^\T H \beta  = \sqrt{n}\beta - \sqrt{n}A^{-1}\left( B + B^\T + C\right) \beta.
		\end{align*}
		Again, by Lemma \ref{lemmaq>1nonbias},
		\begin{align*}
			\left\Vert \sqrt{n} A^{-1}\left( B + B^\T + C\right) \beta \right\Vert \leq \left(n \left\Vert A^{-1} \right\Vert\right) \left(n^{-1/2}\left\Vert B + B^\T + C \right\Vert \right) \left\Vert \beta \right\Vert \conp 0.
		\end{align*}
		Finally, by the Multivariate Central Limit Theorem,
		\begin{align*}
			n^{-1/2} H^\T \epsilon \cond \n_q\left( 0_q , \sigma_\epsilon^2 \Sigma_H \right).
		\end{align*}
		Since $nA^{-1} \conp \Sigma_H^{-1}$, it follows by Slutsky's Theorem that
		\begin{align*}
			\sqrt{n} A^{-1} H^\T \epsilon \cond \n_q\left( 0_q , \sigma_\epsilon^2 \Sigma_H^{-1} \right),
		\end{align*}
		which finishes the proof.
	\end{proof}
	
	\begin{proof}[Proof of Proposition \ref{propositionq>1varianceestimates}]
		This follows from Lemma \ref{lemmaq>1nonbias}.
	\end{proof}

	\subsection{Proofs for Section \ref{sectionsigmamu}}
	
	\begin{proof}[Proof of Proposition \ref{propositionsigmamuefficiency}]
		Letting $\gammahat$ denote the least-squares estimator for $\gamma$, it is known that $\gammahat$ is efficient for estimating $\gamma$ in the low-dimensional linear model.  Since $Z$ is assumed to be of full rank, there exists a smooth re-parameterization of the problem given by $\left( \gamma , \sigma_\epsilon^2 \right) \mapsto \left( \sigmamu, \vartheta, \sigma_\epsilon^2 \right)$, where $\left( \sigmamu,\vartheta\right)$ is the polar representation of $\left\Vert Z_{S_\gamma} \gamma \right\Vert^2$.  Taking the bowl-shaped loss to be quadratic in the first component, the result follows from the arguments of Section 2.3 of \citeA{bickel1993} since $\left\Vert P_{S_\gamma} Y \right\Vert^2 = \left\Vert Z\gammahat \right\Vert^2$.
	\end{proof}
	
	The proof for Theorem \ref{theoremsigmamu} will rely on the proof of Theorem \ref{theoremsigmaepsilon} from Section \ref{sectionproofssigmaepsilon}.
	
	\begin{proof}[Proof of Theorem \ref{theoremsigmamu}]
		Indeed, we may write
		\begin{align*}
			\frac{1}{n} \left\Vert Y \right\Vert^2 = \frac{1}{n} \left\Vert \mu \right\Vert^2 + \frac{2}{n}\mu^\T \epsilon + \frac{1}{n} \left\Vert \epsilon \right\Vert^2
		\end{align*}
		Note that, from equations \eqref{equationsigmaepsilonhats}, \eqref{equationsigmaepsilonhatl} and \eqref{equationsigmaepsilonhatm}, it follows that
		\begin{align*}
			&\sigmamuhats = \frac{1}{n} \left\Vert \mu \right\Vert^2 + \frac{2}{n}\mu^\T \epsilon + \op\left( n^{-1/2} \right),\\
			&\sigmamuhatm = \frac{1}{n} \left\Vert \mu \right\Vert^2 + \frac{2}{n}\mu^\T \epsilon + \op\left( n^{-1/2} \right),\\
			&\sigmamuhatl = \frac{1}{n} \left\Vert \mu \right\Vert^2 + \frac{2}{n}\mu^\T \epsilon + \op\left( n^{-1/2} \right).
		\end{align*}
		By the Multivariate Central Limit Theorem, it follows that
		\begin{align*}
			\sqrt{n}
			\begin{pmatrix} n^{-1} \left\Vert \mu \right\Vert^2 - \sigmamu \\ 2n^{-1} \mu^\T\epsilon \end{pmatrix}
			\cond \n\left( \begin{pmatrix} 0 \\ 0 \end{pmatrix} , \begin{pmatrix} \kappa & 0 \\ 0 & 4\sigma_\epsilon^2 \sigmamu \end{pmatrix} \right).
		\end{align*}
		Applying the Cram\'er-Wold device finishes the proof.
	\end{proof}
	
	\begin{proof}[Proof of Proposition \ref{propositionkappamu}]
		Indeed,
		\begin{align*}
			\kappamuhat =& \frac{1}{n} \sum_{j=1}^n \left( \left( \mu_j^2 - \sigmamu \right) + \left( \muhatj- \mu_j \right)^2 + 2\left( \muhatj - \mu_j \right) \mu_j  - \left( \sigmamuhat - \sigmamu \right) \right)^2\\
			=& \frac{1}{n} \sum_{j=1}^n \left( \mu_j^2 - \sigmamu \right)^2
			+ \frac{1}{n} \sum_{j=1}^n \left( \left( \muhatj - \mu_j \right)^2 + 2\left( \muhatj - \mu_j \right) \mu_j - \left( \sigmamuhat - \sigmamu \right) \right)^2 \\
			&+ \frac{2}{n} \sum_{j=1}^n \left( \mu_j^2 - \sigmamu \right) \left( \left( \muhatj - \mu_j \right)^2 + 2\left( \muhatj - \mu_j \right) \mu_j  - \left( \sigmamuhat - \sigmamu \right) \right).
		\end{align*}
		Applying the Law of Large Numbers yields
		\begin{align}\label{equationkappaconvergence}
			\frac{1}{n} \sum_{j=1}^n \left( \mu_j^2 - \sigmamu \right)^2  \conp \kappa_\mu.
		\end{align}
		By the triangle inequality and Cauchy-Schwarz, it follows that
		\begin{align*}
			\frac{1}{n} \sum_{j=1}^n &\left( \left( \muhatj - \mu_j \right)^2 + 2\left( \muhatj - \mu_j \right) \mu_j - \left( \sigmamuhat - \sigmamu \right) \right)^2 \\
			&\leq \frac{4}{n} \left\Vert \muhat - \mu \right\Vert^4_4
			+ \frac{8}{n} \sum_{j=1}^n \left( \muhatj - \mu_j \right)^2 \mu_j^2
			+ 4 \left( \sigmamuhat - \sigmamu \right)^2\\
			&\leq \frac{4}{n} \left\Vert \muhat - \mu \right\Vert^4_2
			+ \frac{8}{n} \left\Vert \muhat - \mu \right\Vert_2^2 \left\Vert \mu \right\Vert^2_4
			+ 4 \left( \sigmamuhat - \sigmamu \right)^2.
		\end{align*}
		From Theorem \ref{theoremsigmamu}, we see that $\sigmamuhat \conp \sigmamu$.  Therefore, combining this with Proposition \ref{propositionsubgaussianoracle} shows that
		\begin{align}\label{equationsigmamubias}
			\frac{1}{n} \sum_{j=1}^n &\left( \left( \muhatj - \mu_j \right)^2 + 2\left( \muhatj - \mu_j \right) \mu_j - \left( \sigmamuhat - \sigmamu \right) \right)^2 \conp 0.
		\end{align}
		Now, by another application of Cauchy-Schwarz,
		\begin{align*}
			\frac{2}{n} \sum_{j=1}^n &\left|\left( \mu_j^2 - \sigmamu \right) \left( \left( \muhatj - \mu_j \right)^2 + 2\left( \muhatj - \mu_j \right) \mu_j  - \left( \sigmamuhat - \sigmamu \right) \right) \right| \\
			\leq& \frac{2}{n}
			\left( \sum_{j=1}^n \left( \left( \muhatj - \mu_j \right)^2 + 2\left( \muhatj - \mu_j \right) \mu_j - \left( \sigmamuhat - \sigmamu \right) \right)^2 \right)^{1/2}\\
			&\times \left( \sum_{j=1}^n \left( \mu_j^2 - \sigmamu \right)^2 \right)^{1/2} .
		\end{align*}
		From equations \eqref{equationkappaconvergence} and \eqref{equationsigmamubias}, it will follow that
		\begin{align*}
			\frac{2}{n} \sum_{j=1}^n &\left|\left( \mu_j^2 - \sigmamu \right) \left( \left( \muhatj - \mu_j \right)^2 + 2\left( \muhatj - \mu_j \right) \mu_j  - \left( \sigmamuhat - \sigmamu \right) \right) \right|\conp 0.
		\end{align*}
		Combining the results finishes the proof.
	\end{proof}
	
	\subsection{Proofs for Section \ref{sectionsigmaepsilon}}\label{sectionproofssigmaepsilon}
	\begin{proof}[Proof of Theorem \ref{theoremsigmaepsilon}]
		Indeed, note that
		\begin{align}\label{equationsigmaepsilonhats}
			\sigmaepsilonhats
			= \frac{1}{n} \left( \left\Vert \mu - \muhat \right\Vert^2 + \epsilon^\T \left( \mu - \muhat \right) + \left\Vert \epsilon \right\Vert^2 \right)
			= \frac{1}{n} \left\Vert \epsilon \right\Vert^2 + \op\left( n^{-1/2} \right),
		\end{align}
		where the last equality follows from Proposition \ref{propositionsubgaussianoracle} and Lemma \ref{lemmaq=1bias}.  Next, some algebra shows that $\sigmaepsilonhatl$ may be decomposed as
		\begin{align*}
			\sigmaepsilonhatl
			=& \frac{1}{n} \left\Vert \epsilon \right\Vert^2
			+ \frac{1}{n}\sum_{m\in\m_{u_\gamma}} w_{m,Y} \left( \left\Vert \pp{m} \mu \right\Vert^2 + 2\mu^\T\pp{m}\epsilon\right) \\
			&- \frac{1}{n}\sum_{k\in\m_{u_\gamma}} \sum_{m\in\m_{u_\gamma}} w_{k,Y} w_{m,Y} \left( \mu^\T \pp{k} \pp{m} \mu + \epsilon^\T P_k P_m \epsilon \right).			
		\end{align*}
		Applying Corollary \ref{corollarybiasother} yields
		\begin{align*}
			\frac{1}{n}\sum_{m\in\m_{u_\gamma}} w_{m,Y} \left( \left\Vert \pp{m} \mu \right\Vert^2 + 2\mu^\T\pp{m}\epsilon\right) = \op\left( n^{-1/2} \right).
		\end{align*}
		For the other term, it follows from Cauchy-Schwarz, Lemma \ref{lemmabiashw}, and Corollary \ref{corollarybiasother} that
		\begin{align*}
			\frac{1}{n} \sum_{k\in\m_{u_\gamma}} \sum_{m\in\m_{u_\gamma}} w_{k,Y} w_{m,Y} \left( \mu^\T \pp{k} \pp{m} \mu + \epsilon^\T P_k P_m \epsilon \right) = \op\left( n^{-1/2} \right).
		\end{align*}
		Thus, this implies that
		\begin{align}\label{equationsigmaepsilonhatl}
			\sigmaepsilonhatl
			= \frac{1}{n} \left\Vert \epsilon \right\Vert^2 + \op\left( n^{-1/2} \right).
		\end{align}
		Now, by Jensen's inequality,
		\begin{align*}
			\sigmaepsilonhats \leq \sigmaepsilonhatm \leq \sigmaepsilonhatl.
		\end{align*}
		Therefore, we may conclude that
		\begin{align}\label{equationsigmaepsilonhatm}
			\sigmaepsilonhatm = \frac{1}{n} \left\Vert \epsilon \right\Vert^2 + \op\left( n^{-1/2} \right).
		\end{align}
		The asymptotic distribution for all three estimators follows by applying the Central Limit Theorem, which finishes the proof.
	\end{proof}
	
	\begin{proof}[Proof of Corollary \ref{corollarysigmaepsiloncorrelated}]
		Indeed,
		\begin{align*}
			\sigmaepsilonhats = \frac{1}{n} \left\Vert \epsilon \right\Vert^2 + \op\left( n^{-1/2} \right).
		\end{align*}
		From the proof of Lemma \ref{lemmaq=1correlated}, we may apply the Spectral Theorem to obtain the following decomposition: $\Sigma_\epsilon = \Gamma D \Gamma^\T$.  Then, for $\xi \sim \n_n\left( 0_n, I_n \right)$, it follows that $D^{1/2} \xi \equid \epsilon$.  A direct variance calculation for $n^{-1} \left\Vert D^{1/2} \xi \right\Vert^2$ finishes the proof.
	\end{proof}
	
	\begin{proof}[Proof of Proposition \ref{propositionkappaepsilon}]
		The proof is similar to the proof of Proposition \ref{propositionkappamu}.
	\end{proof}

\end{document}